\documentclass[12pt,reqno,oneside]{amsart}
 \usepackage{verbatim,amsmath,amssymb,cite,xspace}
\usepackage{color,cite,graphicx}

 \numberwithin{equation}{section}

\theoremstyle{plain}
\newtheorem{theorem}{Theorem}[section]
\newtheorem{lemma}{Lemma}[section]

\theoremstyle{definition}

\theoremstyle{remark}

\usepackage{color}

\usepackage{amssymb}
\usepackage{amssymb, latexsym}
\usepackage{amsmath}
\usepackage{amscd}
\usepackage{enumerate}
\usepackage{amsfonts}
\usepackage{graphicx}
\usepackage{epsfig}

\newcommand{\R}{\mathbb{R}}

\numberwithin{equation}{section} 

\begin{document}
\title[]{Long time dynamics of defocusing energy critical $3+1$ dimensional wave equation with potential in the radial case }
\author{Hao Jia, Baoping Liu, Guixiang Xu}
\maketitle
\noindent
{\bf Abstract} 
Using channel of energy inequalities developed by T.Duyckaerts, C.Kenig and F.Merle, we prove that, modulo a free radiation, any finite energy radial solution to the defocusing energy critical wave equation with radial potential in $3+1$ dimensions converges to the set of steady states as time goes to infinity. For generic potentials we prove there are only finitely many steady states, and in this case modulo some free radiation the solution converges to one steady state as time goes to infinity.

\begin{section}{Introduction}
 We consider the defocusing energy critical wave equation with a radial potential $V\in C(R^3)$ satisfying $\sup_x(1+|x|)^{\beta}|V(x)|<\infty$ for some $\beta>2$:
\begin{equation}\label{eq:mainequation}
\partial_{tt}u-\Delta u-V(x)u+u^5=0,\,\,{\rm in\,}\,(0,\infty)\times R^3,
\end{equation}
with radial initial data $(u_0,u_1)\in \dot{H}^1\times L^2$. 
Our main result is the following.
\begin{theorem}
Fix $\beta>2$. Let $(u_0,u_1)\in \dot{H}^1\times L^2$ be radial. Define
\begin{equation}
Y:=\{V\in C(R^3):\,V \,\,{\rm is\,\,radial\,\,and\,\,}\sup_x(1+|x|)^{\beta}|V(x)|<\infty\}.
\end{equation}
Denote
\begin{equation}
\Sigma=\{(u_c,0)|\, (u_c,0) \,\,{\rm is\,\,a\,\,radial\,\,steady\,\,state\,\,solution\,\,to\,\,equation\,\,(\ref{eq:mainequation})}\}.
\end{equation}
Let $u\in C([0,\infty),\dot{H}^1)\cap L_t^5L_x^{10}([0,T)\times R^3)$ for any $T<\infty$ be the unique solution to equation (\ref{eq:mainequation}) with initial data $(u(0),\partial_tu(0))=(u_0,u_1)$. Then for some solution $(u^L,\partial_tu^L)$ to linear wave equation without potential we have
\begin{equation}
\lim_{t\to \infty}\inf_{(u_c,0)\in \Sigma}\|(u(t),\partial_tu(t))-(u_c,0)-(u^L(t),\partial_tu^L(t))\|_{\dot{H}^1\times L^2}=0.
\end{equation}
Moreover for $V$ in a dense open set $\Omega\subset Y$, there are only finitely many radial steady states to equation (\ref{eq:mainequation}). In this case, there exist a steady state solution $(u_c,0)$ and some solution $(u^L,\partial_tu^L)$ to linear wave equation without potential, such that
\begin{equation}
\lim_{t\to \infty}\|(u(t),\partial_tu(t))-(u_c,0)-(u^L(t),\partial_tu^L(t))\|_{\dot{H}^1\times L^2}=0.
\end{equation}
\end{theorem}

\medskip
\noindent
Since equation (\ref{eq:mainequation}) is defocusing and energy critical, we know from \cite{Struwe,Gri1,Gri2,BaGe} that the solution $u\in C([0,\infty), \dot{H}^1)\cap L^5_tL^{10}_x([0,T)\times R^3)$ for any $T<\infty$. Moreover the energy
\begin{equation}
E:=\int_{R^3}\frac{|\nabla u|^2}{2}+\frac{(\partial_tu)^2}{2}-\frac{Vu^2}{2}+\frac{u^6}{6}(t,x)\,dx
\end{equation}
is constant along the evolution. Thus the regularity of the equation is well understood. Here the main concern is the dynamics when time goes to infinity. Define the energy functional
\begin{equation}
J(u):=\int_{R^3}\frac{|\nabla u|^2}{2}-\frac{Vu^2}{2}+\frac{u^6}{6}(t,x)\,dx.
\end{equation}
In general when the positive part $V^{+}$ of the potential is large, one can expect that there is a unique positive ground state, which is the global minimizer of energy functional and has negative energy. In addition there can be a number of ``excited states" with higher energies (see Appendix A for more details). It is well known the ground state is asymptotically stable at least when $V$ decays fast. However the dynamics around the excited states can be very complicated even in perturbative regime (and with radial data), involving stable and unstable manifolds. Readers are referred to \cite{TsaiYau,SoffWein,KNS} for a sample of results in this very interesteing area. In particular there might be solutions which stay for very long time near an excited state (just off ``stable manifold") then eventually move away and settle down to either an excited state with lower energy or the ground state. Thus it is extremely hard to follow the solution with all the details. On the other hand, in mathematical physics community, there is widespread belief that solution to dispersive equations should asymptotically decouple into some combination of modulated solitons, free radiation and a term which goes to zero. There are few cases in which this has been proved mathematically, except in perturbative regimes or for integrable equations. Recently a remarkable result was proved in \cite{DKM} for focusing energy critical wave equations in $3+1$ dimensions, which says among many other things that any globally defined radial finite energy solution can be decomposed as the sum of recaled ground state, a free radiation and a term which goes to zero in energy space as time goes to infinity. The proof introduces an important and very natural inequality for radial solution to linear wave equation, which they call ``channel of energy" inequality. More precisely, suppose radial $\overrightarrow{u}^L\in \dot{H}^1\times L^2$ solves the free wave equation
\begin{equation}\label{eq:freewave}
\partial_{tt}u-\Delta u=0,
\end{equation}
with initial data $\overrightarrow{u}^L(0)=(u_0,u_1)$. Then for any $R\ge 0$, the inequality
\begin{equation}\label{eq:energychannel}
\int_{r>R+|t|}^{\infty}(\partial_r(ru^L))^2+(\partial_tru^L)^2(t,r)\,dr\ge \frac{1}{2}\int_R^{\infty}(\partial_r(ru_0))^2+(ru_1)^2\,dr
\end{equation}
holds for all $t\ge 0$ or all $t\leq 0$. We note that (with $S^2$ denoting the $2$ dimensional unit sphere)
\begin{equation}
\int_{|x|\ge R+|t|}|\nabla u^L|^2+(\partial_tu^L)^2(t,x)\,dx\ge|S^2|\int_{r\ge R+|t|}(\partial_r(ru^L))^2+(\partial_tru^L)^2(t,r)\,dr.
\end{equation}
The proof of inequality (\ref{eq:energychannel}) is simple and is based on the observation that $ru$ satisfies $1+1$ dimensional wave equation. Despite its simplicity, this inequality is robust under small perturbations of nonlinearity (the smallness can always be achieved by restriction to the exterior of large light cones), and gives crucial information that is valid for all time. One can expect that near the light cone the dynamics is relatively simple for large times, not only for wave equations, but also for other dispersive equations with finite speed of propagation, such as Klein Gordon equations. Thus there is hope to use the dynamics ``near" lightcone to gain information about the long time dyanamics of solution. The danger is of course that the dynamics near lightcone becomes trivial. In some sense the ``channel of energy" inequality provides the crucial information that rules out this possibility. We note that in the case of Klein Gordon equations where we don't have ``channel of energy" inequalities, similar problem to the one considered here is open. Besides the ``channel of energy" inequality, profile decomposition is another important tool in \cite{DKM}. The ability to write the solution as the sum of large, unit sized and small profiles which evolve almost ``independently" allows the analysis of each profile individually. Since profile decompositions depend heavily on the scale-invariance of the equation and adding potential to the equation destroys such invariance, at first glance it might seem slightly problematic to try to extend the techniques in \cite{DKM} to our setting. However, we observe that outside lightcone, the very small or large profiles are almost unaffected by the potential (see Lemma \ref{lm:profilewithpotential} below). Thus only the unit-sized profile is really influenced by the potential outside any light cone. This enables us to do the usual profile decomposition in the presence of potential outside lightcone. Then we follow closely the arguments in \cite{DKM} to eliminate all the profiles except the ones given by an evolution of free wave and a steady state in the profile decomposition of $u(t_n)$ for any sequence of times $t_n\to \infty$. This is achieved in two main steps. Firstly we show that for any $A>0$ in the region $\{(t,x)\in (0,\infty)\times R^3|\,|x|>t-A\}$ the solution to equation (\ref{eq:mainequation}) behaves almost linearly for large time, given by the free radiation profile. Next, suppose there is another profile $U_{jn}$ besides unit sized profile given by a steady state solution and free wave $\overrightarrow{u}^L(t_n)$. Using the channel of energy inequalities, we show for some positive $\delta$, all $t\ge 0$ or all $t\leq 0$, and $n$ large we have
\begin{equation}
\int_{|x|\ge |t|}|\nabla (u-u^L)|^2+\left(\partial_t(u-u^L)\right)^2(t+t_n,x)dx\ge\frac{1}{2}\int_{|x|\ge |t|}|\nabla U_{jn}|^2+(\partial_tU_{jn})^2(t,x)dx>\delta.
\end{equation}
If this equality holds for all $t\ge 0$, then it contradicts the fact that near lightcone solutions $\overrightarrow{u}(t)$ is almost $\overrightarrow{u}^L(t)$ for large $t$. If this inequality holds for all $t\leq 0$, then we can take $t=-t_n$ and obtain that the initial data $(u_0,u_1)$ has nontrivial energy outside very large ball $B_{t_n}$, again a contradiction. Thus such profile $U_{jn}$ can not exist. Therefore in the profile decomposition of $\overrightarrow{u}(t_n)$ we only have $\overrightarrow{u}^L(t_n)$ and a steady state $(u_c,0)$ as $n\to\infty$ for any sequence of time going to infinity. This settles the proof of the first part of our theorem. For generic potential $V$, we prove in the Appendix there are only finitely many steady states for equation (\ref{eq:mainequation}), thus the second part of the theorem also follows easily. We remark that although it seems unlikely that there would be a continuous family of steady states for equation (\ref{eq:mainequation}), \footnote{We know the radial steady state solutions are uniquely determined by its asymptotic behavior at spatial $\infty$, $u_c\sim\frac{c}{|x|}$. This fact is quite standard, see also \cite{DKM2,KLS} for similar results.} we are not able to rule this possibility out. If this indeed happens for some exotic potential, then modulo free wave some solution might just stay close to the continnum of equilibria without actually converging to any particular one.  \\

\smallskip
\noindent
Comparing with the arguments in \cite{DKM}, we have some slight complications in profile decomposition and ``channel of energy" inequalities for nonlinear solutions due to the presence of potential, which we deal with in Section 2, 4 and the Appendix (which deals with an elliptic problem associated with (\ref{eq:mainequation}) in the similar spirit of \cite{KLS}). On the other hand once we have these tools, the contradiction argument goes much more smoothly, as in the exterior of lightcone, all solutions to equation (\ref{eq:mainequation}) scatter. Thus we don't need to modify certain profiles so that they scatter and painstakingly remove rescaled ground states in the profile decompositions, as has to be done in \cite{DKM}. Such operations are not hard mathematically, however they significantly complicate expositions. Another technical difference is that we use the Strichartz norm $L^5_tL^{10}_x$ to measure scattering, instead of $L^8_{t,x}$\footnote{For the equivalence of these notions of solutions, see \cite{DKM3}.}. The advantage of the space $L^{8}_{t,x}$ is that it's invariant under Lorentz transformations, while $L^5_tL^{10}_x$ is not. However in our problem we don't have Lorentz invariance for equation (\ref{eq:mainequation}). Moreover $u\in L^5_tL^{10}_x$ naturally gives $u^5\in L^1_tL^2_x$, and all the Strichartz estimates follow easily from this without involving any derivatives of $u$. Our paper is organized as follows. In Section 2 we introduce profile decomposition for equation with a potential; in Section 3 we extract linear behavior near lightcone; in Section 4 we recall ``channel of energy" inequalities; in Section 5 we prove the main result; Appendix A and B give necessary facts about the elliptic problem associated with equation (\ref{eq:mainequation}) and the proof of Lemma \ref{lm:othersolutions}.

\end{section}

\begin{section}{some technical lemmas}
In this section we collect some technical lemmas to be used below. We begin with some perturbation results, in the special case when $a\equiv 0$ these results are well known.
\begin{lemma}\label{lm:perturbationresult}
Let $0\in I\subset R$ be an interval of time. Suppose $\tilde{u}(t,x)\in C_t(I,\dot{H}^1(R^3))$ with $\|\tilde{u}\|_{L^5_tL^{10}_x(I\times R^3)}\leq M<\infty$, $\|a\|_{L^{5/4}_tL^{5/2}_x(I\times R^3)}\leq \beta<\infty$ and $e(t,x),\,f(t,x)\in L^1_tL^2_x(I\times R^3)$, satisfy
\begin{equation}
\partial_{tt}\tilde{u}-\Delta \tilde{u}+a(t,x)\tilde{u}+\tilde{u}^5=e,
\end{equation}
with initial data $\overrightarrow{\tilde{u}}(0)=(\tilde{u}_0,\tilde{u}_1)\in \dot{H}^1\times L^2$. Suppose for some sufficiently small positive $\epsilon<\epsilon_0=\epsilon_0(M,\beta)$, 
\begin{equation}
\||e|+|f|\|_{L^1_tL^2_x(I\times R^3)}+\|(u_0,u_1)-(\tilde{u}_0,\tilde{u}_1)\|_{\dot{H}^1\times L^2}<\epsilon.
\end{equation}
Then there is a unique solution $u\in C(I,\dot{H}^1)$ with $\|u\|_{L^5_tL^{10}_x(I\times R^3)}<\infty$, satisfying the equation
\begin{equation}
\partial_{tt}u-\Delta u+a(t,x)u+u^5=f,
\end{equation}
with initial data $\overrightarrow{u}(0)=(u_0,u_1)$. Moreover, we have the following estimate
\begin{equation}
\sup_{t\in I}\|\overrightarrow{u}(t)-\overrightarrow{\tilde{u}}(t)\|_{\dot{H}^1\times L^2}+\|u-\tilde{u}\|_{L^5_tL^{10}_x(I\times R^3)}<C(M,\beta)\epsilon.
\end{equation}
\end{lemma}

\smallskip
\noindent
{\bf Proof.} Pick a small absolute number $\delta$ whose precise value is to be determined. Write $I=\cup_{j=1}^{J}I_j$, where $I_j$ are mutually disjoint and adjacent intervals and $J=J(\beta,M)<\infty$, such that 
\begin{equation}
\sup_j(\|a\|_{L^{5/4}_tL^{5/2}_x(I_j\times R^3)}+\|\tilde{u}\|_{L^5_tL^{10}_x(I_j\times R^3)})<\delta.
\end{equation}
Assume without loss of generality $0\in I_1$. Take a large absolute constant $K>2$ to be determined. Denote $w=u-\tilde{u}$, we shall prove by induction
\begin{equation}\label{eq:induction}
\|w\|_{L^5_tL^{10}_x(I_j\times R^3)}+\sup_{t\in I_j}\|\overrightarrow{w}(t)\|_{\dot{H}^1\times L^2}\leq K^j\epsilon,\,{\rm for\,} 1\leq j\leq J.
\end{equation}
The existence and uniqueness of solution then follow easily from local Cauchy theory and a continuity argument. We shall prove this estimate under the additional assumption that $\|w\|_{L^5_tL^{10}_x(I\times R^3)}\leq K^{J+3}\epsilon$, it is clear how to remove this assumption using continuity argument. By the equations for $u$ and $\tilde{u}$, $w$ must satisfy
\begin{equation}\label{eq:theperturbedequation}
\partial_{tt}w-\Delta w+a(t,x)w+5w\tilde{u}^4+10w^2\tilde{u}^3+10w^3\tilde{u}^2+5w^4\tilde{u}+w^5=-e+f,\,{\rm in\,}I\times R^3,
\end{equation}
with initial data $(w_0,w_1)=\overrightarrow{u}(0)-\overrightarrow{\tilde{u}}(0)$. We can rewrite equation (\ref{eq:theperturbedequation}) as integral equation
\begin{eqnarray*}
&&w(t)=S(t)(w_0,w_1)+\int_0^t\frac{\sin{|\nabla|(t-s)}}{|\nabla|}(-e+f \\
&&\quad\quad\quad\quad\quad-a(x,t)w-5w\tilde{u}^4-10w^2\tilde{u}^3-10w^3\tilde{u}^2-5w^4\tilde{u}-w^5)(s)ds,
\end{eqnarray*}
where $S(t)$ denotes the evolution operator of free wave equation. Thus we obtain by application of Strichartz estimates and H\"{o}lder inequality,
\begin{equation}\label{eq:mainiterationinequality}
\|w\|_{L^5_tL^{10}_x(I_1\times R^3)}+\sup_{t\in I_1}\|w(t)\|_{\dot{H}^1\times L^2}\leq \frac{C_1}{2}\epsilon+\frac{C_1}{2}\epsilon+C_1\left((\delta^4+\delta)\|w\|_{L^5_tL^{10}_x(I_1\times R^3)}+\|w\|_{L^5_tL^{10}_x(I_1\times R^3)}^5\right).
\end{equation}
Here $C_1$ is an absolute constant determined by the constants that appear in Strichartz estimates and H\"{o}lder inequalities. Now first choose an absolute small positive $\delta$ such that
\begin{equation}
C_1(\delta^4+\delta)\leq \frac{1}{4}.
\end{equation}
Fix this $\delta$, then $J=J(M,\beta)$ is determined. Now fix $K=2C_1$, then choose $\epsilon_0=\epsilon_0(M,\beta)>0$ so small that 
\begin{equation}
C_1(K^{J+3}\epsilon_0)^4\leq \frac{1}{4}.
\end{equation}
Then we have by the inequality (\ref{eq:mainiterationinequality}), 
\begin{equation}
\|w\|_{L^5_tL^{10}_x(I_1\times R^3)}+\sup_{t\in I_1}\|w(t)\|_{\dot{H}^1\times L^2}\leq 2C_1\epsilon\leq K\epsilon.
\end{equation}
Suppose the estimate $\|w\|_{L^5_tL^{10}_x(I_j\times R^3)}+\sup_{t\in I_j}\|\overrightarrow{w}(\cdot,t)\|_{\dot{H}^1\times L^2}\leq K^j\epsilon$ is true for some $j\ge 1$, we shall prove it remains true for $j+1$. Again by writing equation (\ref{eq:theperturbedequation}) as integral equation starting from $t_j$ which separates $I_j$ and $I_{j+1}$, using the estimates on $I_j$ (actually we only need the estimate at the point $t_j$ separating $I_j$ and $I_{j+1}$), we obtain
\begin{eqnarray*}
&&\|w\|_{L^5_tL^{10}_x(I_{j+1}\times R^3)}+\sup_{t\in I_{j+1}}\|w(t)\|_{\dot{H}^1\times L^2}\\
&&\quad\leq C_1K^j\epsilon+C_1\left((\delta^4+\delta)\|w\|_{L^5_tL^{10}_x(I_{j+1}\times R^3)}+\|w\|_{L^5_tL^{10}_x(I_{j+1}\times R^3)}^5\right).
\end{eqnarray*}
Again by the choice of parameters, we obtain
\begin{eqnarray*}
&&\|w\|_{L^5_tL^{10}_x(I_{j+1}\times R^3)}+\sup_{t\in I_{j+1}}\|w(t)\|_{\dot{H}^1\times L^2}\\
&&\quad\leq \frac{C_1}{2}K^j\epsilon+\frac{C_1}{2}\epsilon+\frac{1}{2}\|w\|_{L^5_tL^{10}_x(I_{j+1}\times R^3)}+\frac{1}{2}\sup_{t\in I_{j+1}}\|w(t)\|_{\dot{H}^1\times L^2}.
\end{eqnarray*}
Thus $\|w\|_{L^5_tL^{10}_x(I_{j+1}\times R^3)}+\sup_{t\in I_{j+1}}\|w(t)\|_{\dot{H}^1\times L^2}\leq 2C_1K^j\epsilon\leq K^{j+1}\epsilon$. The proof is finished.

\begin{lemma}\label{lm:globalregularityofperturbedwave}
Let $I$ be an interval of time and $a\in L^{5/4}_tL^{5/2}_x\cap L^1_tL^3_x(I\times R^3)$, and $f\in L^1_tL^2_x(I\times R^3)$, with bounds $\|a\|_{L^{5/4}_tL^{5/2}_x}+\|a\|_{L^1_tL^3_x}\leq M$ and $\|f\|_{L^1_tL^2_x}\leq \beta$. Then there exists a unique solution $u\in C(I,\dot{H}^1)\cap L^5_tL^{10}_x(I\times R^3)$ to the equation
\begin{equation}\label{eq:perturbedwave}
\partial_{tt}u-\Delta u+a(t,x)u+u^5=f,
\end{equation} 
with initial data $(u_0,u_1)\in \dot{H}^1\times L^2$ ($\|(u_0,u_1)\|_{\dot{H}^1\times L^2}\leq E$). Moreover, we have
\begin{equation}\label{eq:bound}
\|u\|_{L^5_tL^{10}_x}(I\times R^3)\leq C(E,M,\beta).
\end{equation}
Thus if $I=R$, then there exist solutions $u^L_{+},\,u^L_{-}$ to free wave equation, such that
\begin{eqnarray}
&&\lim_{t\to +\infty}\|u(t)-u^L_{+}(t)\|_{\dot{H}^1\times L^2}=0,\\
&&\lim_{t\to -\infty}\|u(t)-u^L_{-}(t)\|_{\dot{H}^1\times L^2}=0.
\end{eqnarray}
\end{lemma}

\smallskip
\noindent
{\bf Proof.} We only need to prove the bound (\ref{eq:bound}), the other statements follow easily from well known arguments. Again we prove this estimate as a priori estimate, assuming $\|u\|_{L^5_tL^{10}_x}(I\times R^3)<\infty$ and $u$ is smooth with fast decay. Then it's easy to remove this assumption by standard local Cauchy theory, approximation and continuity argument. We first study the energy inequality. Multiplying $\partial_tu$ to equation (\ref{eq:perturbedwave}) and integrate, we obtain
\begin{eqnarray*}
&&\frac{d}{dt}\int_{R^3}\frac{(\partial_tu)^2}{2}+\frac{|\nabla u|^2}{2}+\frac{u^6}{6}(t,x)dx\\
&&\quad=\int_{R^3}f\partial_tu(t,x)dx-\int_{R^3}au\partial_tu(t,x)dx\\
&&\quad\leq C\left(\int_{R^3}(\partial_tu)^2(t,x)dx\right)^{\frac{1}{2}}\left(\int_{R^3}f^2(t,x)dx\right)^{\frac{1}{2}}\\
&&\quad\quad+C\left(\int_{R^3}|a|^3(t,x)dx\right)^{\frac{1}{3}}\left(\int_{R^3}u^6(t,x)dx\right)^{\frac{1}{6}}\left(\int_{R^3}(\partial_tu)^2(t,x)dx\right)^{\frac{1}{2}}\\
&&\quad\leq C\left(\int_{R^3}(\partial_tu)^2(t,x)dx\right)\|f(t)\|_{L^2_x}+C\|f(t)\|_{L^2_x}\\
&&\quad\quad+C\|a(t)\|_{L^3_x}\left(\int_{R^3}|\nabla u|^2(t,x)dx\right)+C\|a(t)\|_{L^3_x}\left(\int_{R^3}(\partial_tu)^2(t,x)dx\right)
\end{eqnarray*}
Since $f\in L^1_tL^2_x(I\times R^3),\,a\in L^1_tL^3_x(I\times R^3)$, we can use Gronwall's inequality and conclude
\begin{equation}
\sup_{t\in I}\|\overrightarrow{u}(t)\|_{\dot{H}^1\times L^2}\leq C(E,M,\beta).
\end{equation}
Write $I=\cup_{j=1}^JI_j$ where $I_j$ are adjacent and mutually disjoint intervals, such that
\begin{equation}
\|a\|_{L^{5/4}_tL^{5/2}_x(I_j\times R^3)}+\|f\|_{L^1_tL^2_x(I_j\times R^3)}\leq \delta\,\,{\rm for\,each\,}j,
\end{equation}
where $\delta=\delta(M,E,\beta)>0$ is to be determined. Thus $J=J(M,E,\beta)$. Take $t_j\in I_j$, and solution $\tilde{u}_j$ to
\begin{equation}
\partial_{tt}\tilde{u}_j-\Delta \tilde{u}_j+\tilde{u}_j^5=0\,\,{\rm in\,}I_j\times R^3,
\end{equation}
with $\overrightarrow{\tilde{u}_j}(t_j)=\overrightarrow{u}(t_j)$. By Corollary 2 in \cite{BaGe}, we have 
\begin{equation}
\|\tilde{u}_j\|_{L^5_tL^{10}_x(I_j\times R^3)}\leq C(E,M,\beta).
\end{equation}
See also in \cite{Tao} for an explicit bound for the spacetime norm. Thus in $I_j\times R^3$, $\tilde{u}_j$ satisfies
\begin{equation}
\partial_{tt}\tilde{u}_j-\Delta \tilde{u}_j+a(t,x)\tilde{u}_j+\tilde{u}_j^5=a(t,x)\tilde{u}_j\,\,{\rm in\,}I_j\times R^3,
\end{equation}
with $\overrightarrow{\tilde{u}_j}(t_j)=\overrightarrow{u}(t_j)$. Now choose $\delta>0$ sufficiently small so that
\begin{eqnarray*}
\|a\tilde{u}_j\|_{L^1_tL^2_x(I_j\times R^3)}+\|f\|_{L^1_tL^2_x(I_J\times R^3)}\leq \delta C(E,M,\beta)<\epsilon_0(E,M,\beta),
\end{eqnarray*}
where $\epsilon_0$ is given by Lemma \ref{lm:perturbationresult}. Thus we can apply Lemma \ref{lm:perturbationresult} and conclude on each $I_j$, 
\begin{equation}
\|u\|_{L^5_tL^{10}_x(I_j\times R^3)}\leq C(E,M,\beta).
\end{equation}
Sum over $j$, the lemma is proved.\\

\smallskip
\noindent
{\bf Remarks.} In applications below to equation (\ref{eq:main}) $V$ does not decay in time. However if we only consider property of solutions in the region $S=\{(t,x):|x|>|t|\}$, noting that
$\|V\|_{L^{5/4}_tL^{5/2}_x\cap L^1_tL^3_x(S)}<\infty$ and finite speed of propagation, we can apply Lemma \ref{lm:globalregularityofperturbedwave} in $S$.\\

\smallskip
\noindent
Next we make precise the statement that for initial data with large profiles the influence of potential is small.
\begin{lemma}\label{lm:aux}
Let $a\in L^{5/4}_tL^{5/2}_x(R\times R^3)$ and $U^L$ be a solution to the free wave equation in $R\times R^3$. Take parameters $(\lambda_n,t_n)$ with $\lambda_n>0,\,t_n\in R$. Assume one of the following conditions holds:\\
1. $t_n\equiv 0$, $\lim\limits_{n\to\infty}(\lambda_n+\frac{1}{\lambda_n})=\infty$,\\
2. $\lim\limits_{n\to\infty}\frac{t_n}{\lambda_n}\in\{\pm\infty\}$.\\
Let $U$ be the nonlinear profile associated with $U^L,\,\lambda_n,\,t_n$. More precisely
\begin{equation}
\partial_{tt}U-\Delta U+U^5=0\,\,{\rm in\,}R\times R^3,
\end{equation}
with $\overrightarrow{U}(0)=(U^L(0),\partial_tU^L(0))$\, if\, $t_n\equiv 0$; or with 
\begin{equation}
\lim_{t\to +\infty}\|\overrightarrow{U}(t)-\overrightarrow{U^L}(t)\|_{\dot{H}^1\times L^2}=0,\,\,(\lim_{t\to-\infty})
\end{equation}
 if\, $\lim_{n\to\infty}\frac{t_n}{\lambda_n}=-\infty$ ($\lim=\infty$ respectively). Let $u_n$ be the solution to the Cauchy problem
\begin{equation}
\partial_{tt}u_n-\Delta u_n+a(t,x)u_n+u_n^5=0\,\,{\rm in\,}R^3\times R,
\end{equation}
with $\overrightarrow{u}_n(0)=\left(\frac{1}{\lambda_n^{1/2}}U^L(-\frac{t_n}{\lambda_n},\frac{x}{\lambda_n}),\frac{1}{\lambda_n^{3/2}}\partial_tU^L(-\frac{t_n}{\lambda_n},\frac{x}{\lambda_n})\right)$. Then
\begin{equation}
\lim_{n\to\infty}\left(\sup_{t\in R}\|\overrightarrow{u_n}(t)-\overrightarrow{U_n}(t)\|_{\dot{H}^1\times L^2}+\|u_n-U_n\|_{L^5_tL^{10}_x(R\times R^3)}\right)=0,
\end{equation}
where $\overrightarrow{U}_n(x,t)=(\frac{1}{\lambda_n^{1/2}}U(\frac{t-t_n}{\lambda_n},\frac{x}{\lambda_n}),\frac{1}{\lambda_n^{3/2}}\partial_tU(\frac{t-t_n}{\lambda_n},\frac{x}{\lambda_n}))$.
\end{lemma}
\smallskip
\noindent
{\bf Proof.} By the definition of $U_n$, it satisfies
\begin{equation}
\partial_{tt}U_n-\Delta U_n+a(t,x)U_n+U_n^5=a(t,x)U_n\,\,{\rm in}\,\,R\times R^3,
\end{equation}
and 
\begin{equation}
\lim_{n\to \infty}\left\|\overrightarrow{U}_n(0)-\left(\frac{1}{\lambda_n^{1/2}}U^L(-\frac{t_n}{\lambda_n},\frac{x}{\lambda_n}),\frac{1}{\lambda_n^{3/2}}\partial_tU^L(-\frac{t_n}{\lambda_n},\frac{x}{\lambda_n})\right)\right\|_{\dot{H}^1\times L^2}=0.
\end{equation}
Since $U\in L^5_tL^{10}_x(R^3\times R)$ and $U_n(t,x)=\frac{1}{\lambda_n^{1/2}}U(\frac{t-t_n}{\lambda_n},\frac{x}{\lambda_n})$, by the condition 1 or 2 on the scaling and time translation parameters, it's straightforward to verify that
\begin{equation}
\lim_{n\to\infty}\|aU_n\|_{L^1_tL^2_x(R\times R^3)}=0.
\end{equation}
Thus this lemma follows directly from Lemma \ref{lm:perturbationresult}.\\

We now introduce the profile decompositions for wave equation with potential. 
\begin{lemma}\label{lm:profilewithpotential}
Let $a\in L^{5/4}_tL^{5/2}_x\cap L^{1}_tL^3_x(R\times R^3)$. Suppose  radial $(u_{0n},u_{1n})\in \dot{H}^1\times L^2$ are uniformly bounded and have the following linear profile decompositions
\begin{equation}
(u_{0n},u_{1n})=\overrightarrow{U}^L_1(0)+\sum_{j=2}^J(\frac{1}{\lambda_{jn}^{1/2}}U^L_j(-\frac{t_{jn}}{\lambda_{jn}},\frac{x}{\lambda_{jn}}),\frac{1}{\lambda_{jn}^{3/2}}\partial_tU^L_j(-\frac{t_{jn}}{\lambda_{jn}},\frac{x}{\lambda_{jn}}))+\overrightarrow{w}_{Jn}(0),
\end{equation}
with the following properties:\\
\begin{eqnarray*}
&&U^L_j,\,w_{Jn} \,{\rm \,are\,\,radial\,\,and\, \,solve\,\,the\,\, free\, \,wave\,\, equation}\,\,{\rm for\,\,each}\,\,j,\,J;\\
&&{\rm either}\,\,t_{jn}\in R,\,\lambda_{jn}>0,\,\lim_{n\to\infty}\frac{t_{jn}}{\lambda_{jn}}\in\{\pm\infty\}{\rm \,\,or\,\,}t_{jn}\equiv 0,\,\lim_{n\to\infty}\left(\lambda_{jn}+\frac{1}{\lambda_{jn}}\right)=\infty;\\
&& {\rm for\,}\,j\neq j',\,\lim_{n\to\infty}\left(\frac{\lambda_{jn}}{\lambda_{j'n}}+\frac{\lambda_{j'n}}{\lambda_{jn}}+\frac{|t_{jn}-t_{j'n}|}{\lambda_{jn}}\right)=\infty;\\
&&{write}\,\,w_{Jn}(t,x)=\frac{1}{\lambda_{jn}^{1/2}}\tilde{w}^j_{Jn}(\frac{t-t_{jn}}{\lambda_{jn}},\frac{x}{\lambda_{jn}}),\,\,{then}\,\,\tilde{w}_{Jn}^j\rightharpoonup 0,\,\,{\rm and\,}\,w_{Jn}\rightharpoonup 0,\,\,{\rm as}\,\,n\to\infty;\\
&&\lim_{J\to\infty}\limsup_{n\to\infty}\|w_{Jn}\|_{L^5_tL^{10}_x(R^3\times R)}=0.
\end{eqnarray*}
Let $U_1$ satisfy
\begin{equation}
\partial_{tt}U_1-\Delta U_1+a(t,x)U_1+U_1^5=0,{\rm \,in\,}R\times R^3,
\end{equation}
with $\overrightarrow{U}_1(0)=\overrightarrow{U}^L_1(0)$.
Let $U_j$ be the nonlinear profile associated to $U^L_j,\,\lambda_{jn},\,t_{jn}$ as defined in Lemma \ref{lm:aux} for $j\ge 2$. Let $u_n$ be the solution to
\begin{equation}
\partial_{tt}u_n-\Delta u_n+a(t,x)u_n+u_n^5=0,\,{\rm in\,}R\times R^3,
\end{equation}
with $\overrightarrow{u}_n(0)=(u_{0n},u_{1n})$. Then we have the following decomposition 
\begin{equation}\label{eq:decomposition}
\overrightarrow{u}_n(t)=\overrightarrow{U_1}(t)+\sum_{j=2}^J \overrightarrow{U}_{jn}(t)+\overrightarrow{w}_{Jn}(t)+\overrightarrow{r}_{Jn}(t),
\end{equation}
with 
\begin{equation}\label{eq:errorterm}
\lim_{J\to\infty}\limsup_{n\to\infty}\left(\sup_{t\in R}\|\overrightarrow{r}_{Jn}(t)\|_{\dot{H}^1\times L^2}+\|r_{Jn}\|_{L^5_tL^{10}_x(R\times R^3)}\right)=0,
\end{equation}
where $\overrightarrow{U}_{jn}(t,x)=(\frac{1}{\lambda_{jn}^{1/2}}U_j(\frac{t-t_{jn}}{\lambda_{jn}},\frac{x}{\lambda_{jn}}),\frac{1}{\lambda_{jn}^{3/2}}\partial_tU_j(\frac{t-t_{jn}}{\lambda_{jn}},\frac{x}{\lambda_{jn}}))$. Moreover, denoting $U_{1n}=U_1$, for $\rho_n>\sigma_n>0$ and $\theta_n\in R$ we have the following orthogonality property for $1\leq j\not= j'$
\begin{eqnarray}
&&\label{eq:no1}\lim_{n\to\infty}\int_{\sigma_n<|x|<\rho_n}\nabla U_{jn}\nabla U_{j'n}+\partial_tU_{jn}\partial_tU_{j'n}(\theta_n,x)dx=0;\\
&&\label{eq:no2}\lim_{n\to\infty}\int_{\sigma_n<|x|<\rho_n}\nabla U_{jn}\nabla w_{Jn}+\partial_tU_{jn}\partial_tw_{Jn}(\theta_n,x)dx=0.
\end{eqnarray}
\end{lemma}

\smallskip
\noindent
{\bf Proof.}  We only need to prove the error estimate (\ref{eq:errorterm}). Once this is done, the other claims follow from similar arguments as in the Appendix B of \cite{DKM} (see also remarks below). Let
\begin{equation}
\overrightarrow{u}_{Jn}:=\overrightarrow{U_1}+\sum_{j=2}^J \overrightarrow{U}_{jn}+\overrightarrow{w}_{Jn}.
\end{equation}
Then $u_{Jn}$ satisfies
\begin{equation}
\partial_{tt}u_{Jn}-\Delta u_{Jn}+a(t,x)u_{Jn}+u_{Jn}^5=f_{Jn}, \,\,{\rm in\,}\,R\times R^3,
\end{equation}
where 
\begin{equation}
f_{Jn}=\sum_{j=2}^Ja(t,x)U_{jn}+a(t,x)w_{Jn}+u_{Jn}^5-\sum_{j=2}^JU_{jn}^5-U_1^5.
\end{equation}
By the definition of nonlinear profile $U_j$, it's clear
\begin{equation}
\lim_{n\to\infty}\|\overrightarrow{u}_{Jn}(0)-\overrightarrow{u}_n(0)\|_{\dot{H}^1\times L^2}=0.
\end{equation}
To prove the error estimate on $\overrightarrow{u}_{Jn}(t)-\overrightarrow{u}_n(t)$, by Lemma (\ref{lm:perturbationresult}), it suffices to show
\begin{equation}
\lim_{J\to\infty}\limsup_{n\to\infty}\|f_{Jn}\|_{L^1_tL^2_x(R\times R^3)}=0.
\end{equation}
This property follows immediately from the assumptions on the parameters $\lambda_{jn},\,t_{jn}$ and on $w_{Jn}$.\\

\smallskip
\noindent
{\bf Remarks.} It's easy to see from equations (\ref{eq:no1},\ref{eq:no2}) that
\begin{eqnarray}
&&\lim_{n\to\infty}\sup_{t\in R}\left|\int_{\sigma_n<|x|<\rho_n}\nabla U_{jn}\nabla U_{j'n}+\partial_tU_{jn}\partial_tU_{j'n}(t,x)dx\right|=0;\\
&&\lim_{n\to\infty}\sup_{t\in R}\left|\int_{\sigma_n<|x|<\rho_n}\nabla U_{jn}\nabla w_{Jn}+\partial_tU_{jn}\partial_tw_{Jn}(t,x)dx\right|=0.
\end{eqnarray}
We can briefly recall the ideas of proof of orthogonality property of profiles. For profiles $U_{jn}$ and $U_{j'n}$, if $\lambda_{jn}$ and $\lambda_{j'n}$ are not comparable as 
$n\to\infty$ then the orthogonality property follows immediately as the two profiles are supported at different \emph{scales}. If $\lambda_{jn}$ and $\lambda_{j'n}$ are comparable, then by the assumption on the parameters $\lim\limits_{n\to\infty}\frac{|t_{jn}-t_{j'n}|}{\lambda_{jn}}=\infty$. A moment of reflection involving the support properties of profiles then shows the only times the two profiles have nontrivial overlap are when $|\theta_n- \frac{t_{jn}+t_{j'n}}{2}|= O(\lambda_{jn})$. By this $\theta_n$, the two profiles have already evolved for very long time and can be treated as linear profile by scattering results. Thus we are reduced to consider linear profiles with the same parameters only. The orthogonality for $U_{jn}$ and $w_{Jn}$ is slightly different, as there is no characteristic length scale for $w_{Jn}$. One can think as follows. Firstly by appropriate rescaling and time translations, it suffices to consider the case of unit sized profile with $t_{jn}\equiv 0$. Then by our assumption (appropriately rescaled and translated) $w_{Jn}\rightharpoonup 0$ as $n\to\infty$, it suffices to consider the case $\theta_n\to\infty$. In this case again the profile has evoloved for very long time by $\theta_n$ so that we can use a linear profile to replace it. Thus in summary it suffices to consider linear profiles $\tilde{U}^L_{jn}$ and $w_{Jn}$. Since we are in radial case the proof then follows from relatively straightforward calculations. We only note that in the calculations it is helpful to first rescale and translate one family of profiles in time so that they are of unit size and starting from time $0$ (the other family then weakly goes to zero).
\end{section}

\begin{section}{Existence of free wave}
Now we begin to study defocusing energy critical wave equation with radial potential $V\in Y$:
\begin{equation}\label{eq:main}
\partial_{tt}u-\Delta u-Vu+u^5=0,\,\,{\rm in}\,\,R\times R^3,
\end{equation}
with radial initial data $(u_0,u_1)\in \dot{H}^1\times L^2$. It is now well known that equation (\ref{eq:main}) is globally wellposed. Moreover,
the energy
\begin{equation}
E:=\int_{R^3}\frac{(\partial_tu)^2}{2}+\frac{|\nabla u|^2}{2}-\frac{Vu^2}{2}+\frac{u^6}{6}(t,x)dx
\end{equation}
is constant for $t\in R$. Thus 
\begin{equation}
\sup_{t\in R}\|\overrightarrow{u}(t)\|_{\dot{H}^1\times L^2}\leq C(E,\|V\|_Y).
\end{equation}
In this section, following arguments in \cite{DKM} we show for large time near light cone, $u(t)$ behaves almost linearly.

\begin{lemma}\label{lm:inter}
Suppose $u$ is the solution to equation (\ref{eq:main}) with radial finite energy initial data $(u_0,u_1)$. For any $A>0$, there exists a radial solution $v^L_A$ to the free wave equation, such that
\begin{equation}
\lim_{t\to\infty}\|\overrightarrow{u}(t)-\overrightarrow{v}^L_A(t)\|_{\dot{H}^1\times L^2(|x|\ge t-A)}=0.
\end{equation}
\end{lemma}

\smallskip
\noindent
{\bf Proof.} Take a sequence of time $t_n\to\infty$. Define for $(t,x)\in R\times R^3$
\begin{eqnarray}
V_n(t,x)=\left\{\begin{array}{rl}
        V(|x|)& {\rm if}\,\,|x|\ge |t-t_n|;\\
        V(|t-t_n|)& {\rm if}\,\,|x|<|t-t_n|.
        \end{array}\right.
\end{eqnarray}
Let $u_n$ be solution to 
\begin{equation}
\partial_{tt}\phi-\Delta\phi+V_n\phi+\phi^5=0, \,\,{\rm in}\,\,(t_n,\infty)\times R^3,
\end{equation}
with $\overrightarrow{u}_n(t_n)=(u(t_n),\partial_tu(t_n))$. By Lemma \ref{lm:globalregularityofperturbedwave} we know that $u_n$ exists globally and scatters. Thus we can find solution $u^L_n$ 
to the free wave equation such that
\begin{equation}
\lim_{t\to\infty}\|\overrightarrow{u}_n(t)-\overrightarrow{u}^L_n(t)\|_{\dot{H}^1\times L^2}=0.
\end{equation}
By finite speed of propagation property of wave equation and the definition of $V_n$, we see 
\begin{equation}
\overrightarrow{u}(t,x)=\overrightarrow{u_n}(t,x)\,\,{\rm for}\,\,|x|\ge t-t_n {\rm \,\,when\,}\,t\ge t_n.
\end{equation}
Thus
\begin{equation}
\lim_{t\to\infty}\|\overrightarrow{u}(t)-\overrightarrow{u}^L_n(t)\|_{\dot{H}^1\times L^2(|x|\ge t-t_n)}=0.
\end{equation}
Take $t_N$ with $t_N\ge A$, we see the lemma is proved by setting $\overrightarrow{v}^L_A=\overrightarrow{u}^L_N$.

\begin{theorem}\label{th:freewave}
Let $u,\,V$ be as in Lemma \ref{lm:inter}. Then there exists a unique solution $U^L$ to the free wave equation which is radial, such that
\begin{equation}
\lim_{t\to\infty}\|\overrightarrow{u}(t)-\overrightarrow{U}^L(t)\|_{\dot{H}^1\times L^2(|x|\ge t-A)}=0,\,\,{\rm for\,\,any}\,\,A\ge 0.
\end{equation}
\end{theorem}

\smallskip
\noindent
{\bf Proof.} Take a sequence of time $t_n\to\infty$ and let $(u_{0n},u_{1n})=(u(t_n),\partial_tu(t_n))$. Passing to a subsequence if necessary, we may assume $(u_{0n},u_{1n})$
has the following profile decomposition
\begin{equation}
(u_{0n},u_{1n})=(U^L(t_n),\partial_tU^L(t_n))+\sum_{j=2}^J\overrightarrow{U}^L_{j,n}(0)+\overrightarrow{w}_{Jn}(0),
\end{equation}
where 
\begin{equation}
\overrightarrow{U}^L_{jn}=(\frac{1}{\lambda_{jn}^{1/2}}U^L_j(-\frac{t_{jn}}{\lambda_{jn}},\frac{x}{\lambda_{jn}}),\frac{1}{\lambda_{jn}^{3/2}}\partial_tU^L_j(-\frac{t_{jn}}{\lambda_{jn}},\frac{x}{\lambda_{jn}})),
\end{equation}
 and the parameters $(\lambda_{jn},\,t_{jn})$ and $(1,t_n)$ satisfy the usual orthogonality properties. For any $A>0$, let $u^L_A$ be a solution to the free wave equation as given by Lemma \ref{lm:inter}. 
We then have the following profile decompositions for $(u_{0n},u_{1n})-(u_A^L(t_n),\partial_t u_A^L(t_n))$
\begin{equation}
(u_{0n},u_{1n})-(u^L_A(t_n),\partial_t u^L_A(t_n))=(U^L(t_n)-u^L_A(t_n),\partial_tU^L(t_n)-\partial_tu^L_A(t_n))+\sum_{j=2}^J\overrightarrow{U}^L_{j,n}(0)+\overrightarrow{w}_{Jn}(0).
\end{equation}
By the orthogonality property of profiles (see \cite{DKM} and the remarks below Lemma \ref{lm:profilewithpotential}) we obtain
\begin{eqnarray*}
&&\limsup_{n\to\infty}\|(U^L(t_n)-u^L_A(t_n),\partial_tU^L(t_n)-\partial_tu^L_A(t_n))\|_{\dot{H}^1\times L^2(|x|\ge t_n-A)}\\
&&\quad\leq\limsup_{n\to\infty}\|(u_{0n},u_{1n})-(u^L_A(t_n),\partial_t u^L_A(t_n))\|_{\dot{H}^1\times L^2(|x|\ge t_n-A)}\\
&&{\rm thus\,\,by\,\,definition\,\,of\,\,}u^L_A\,\,\\
&&\quad\leq\limsup_{n\to\infty}\|(u(t_n),\partial_tu(t_n))-(u^L_A(t_n),\partial u^L_A(t_n))\|_{\dot{H}^1\times L^2(|x|\ge t_n-A)}\\
&&\quad = 0.
\end{eqnarray*}
Since $\overrightarrow{U}^L-\overrightarrow{u}^L_A$ is solution to free wave equation, the quantity
\begin{equation}
\|\overrightarrow{U}^L(t)-\overrightarrow{u}^L_A(t)\|_{\dot{H}^1\times L^2(|x|\ge t-A)} 
\end{equation}
is nonincreasing in $t$, we obtain
\begin{equation}
\lim_{t\to\infty}\|\overrightarrow{U}^L(t)-\overrightarrow{u}^L_A(t)\|_{\dot{H}^1\times L^2(|x|\ge t-A)} =0.
\end{equation}
Thus
\begin{equation}
\lim_{t\to\infty}\|\overrightarrow{u}(t)-\overrightarrow{U}^L(t)\|_{\dot{H}^1\times L^2(|x|\ge t-A)} =0.
\end{equation}
Since $A>0$ is arbitrary, the existence part of the theorem is proved. Uniqueness follows from strong Huygens principle.
\end{section}

\begin{section}{channel of energy inequalities}
We shall use the following channel of energy inequalities (see \cite{DKM1} and Appendix C of \cite{DKM} for proofs).
\begin{lemma}\label{lm:channelofenergyforlinearwave}
Suppose $u^L$ satisfies the free wave equation with radial initial data $(u_0,u_1)\in \dot{H}^1\times L^2$. Then for any $r_0\ge 0$, the following ``channel of energy" inequality
holds for all $t\ge 0$ or all $t\leq 0$:
\begin{equation}
\int_{r\ge r_0+|t|}\left(\partial_r(ru^L)\right)^2(t,r)+(\partial_t ru^L)^2(t,r)dr\ge \frac{1}{2}\int_{r\ge r_0}\left(\partial_r(ru_0)\right)^2+u_1^2dr.
\end{equation}
\end{lemma}

\smallskip
\noindent
{\bf Remarks.} We refer readers to \cite{DKM} for details. We only remark the proof is based on the observation that $ru^L(t,r)$ solves $1+1$ dimensional wave equation and thus
\begin{equation}
ru^L(t,r)=f(r-t)+g(r+t),
\end{equation}
for some $f$ and $g$. The rest is mostly direct calculations.

\begin{lemma}\label{lm:channelwithdirection}
Suppose $u^L$ satisfies the free wave equation with radial initial data $(u_0,u_1)\in \dot{H}^1\times L^2$. Then there exists $t_0\ge 0$, such that for $t\ge t_0$, we have
\begin{equation}
\int_{r\ge t-t_0}\left(\partial_r(ru^L)\right)^2(t,r)+(\partial_t ru^L)^2(t,r)dr\ge \frac{1}{2}\int_{r>0}\left(\partial_r(ru^L)\right)^2(t_0,r)+(\partial_t ru^L)^2(t_0,r)dr.
\end{equation}
\end{lemma}

\smallskip
\noindent
{\bf Remarks.} Again we refer readers to \cite{DKM} for details. We only remark that since we can extend $ru^L(t,|r|)$ as an odd function in $r\in R$, we have the following formula
\begin{equation} 
ru^L(t,|r|)=f(r-t)-f(-r-t),
\end{equation}
for some $f$. The rest is mostly direct calculations.\\

\smallskip
\noindent
We shall also need the following lemma on the ``growth" of support for solutions to linear wave equation with potential in the radial case.
\begin{lemma}\label{lm:propagationofsupport}
Suppose $a\in L^{5/4}_tL^{5/2}_x(K)$ for any $K\Subset R\times R^3$. Let $u$ be the solution to
\begin{equation}
\partial_{tt}u-\Delta u+a(t,x)u=0,
\end{equation}
with compactly supported radial initial data $(u_0,u_1)\in\dot{H}^1\times L^2$. Denote
\begin{equation}
\rho(f,g):=\inf\{r>0| \,{\rm supp\,}(f,g)\subseteq B(0,r)\}.
\end{equation}
Then for all $t\ge 0$ or all $t\leq 0$, we have
\begin{equation}
\rho(u(t),\partial_tu(t))=\rho(u_0,u_1)+|t|.
\end{equation}
\end{lemma}

\smallskip
\noindent
{\bf Proof.} The proof is almost the same as the proof of part (a) of Proposition 2.2 in \cite{DKM}, we only sketch some of the details below. Let $u^L$ be the solution to free wave equation with initial data $\overrightarrow{u}(0)=(u_0,u_1)$. By the ``channel of energy" inequality for $u^L$, without loss of generality, we consider the case that for $t\ge 0$, there is a sequence $\rho_n\to\rho(u_0,u_1)-$, such that
\begin{eqnarray*}
&&\frac{1}{|S^2|}\int_{|x|\ge \rho_n+t}\left(\frac{|\nabla u^L|^2}{2}+\frac{(\partial_tu^L)^2}{2}\right)(t,x)dx\\
&&\quad\ge\int_{r\ge \rho_n+t}\left(\partial_r(ru^L)\right)^2+(\partial_tru^L)^2(t,r)dr\\
&&\quad{\rm by \,\, channel\,\,of\,\,energy\,\,inequality\,\,}\\
&&\quad\ge \frac{1}{2}\int_{r\ge \rho_n}\left(\partial_r(ru_0)\right)^2+(ru_1)^2(r)dr\\
&&\quad\ge \frac{1}{2}\left(\int_{r\ge \rho_n}\left(r\partial_ru_0\right)^2+(ru_1)^2(r)dr-\rho_nu_0(\rho_n)^2\right)\\
&&\quad\ge\frac{1}{2}\left(\frac{1}{|S^2|}\int_{|x|\ge \rho_n}|\nabla u_0|^2+u_1^2dx-\rho_n\left(\frac{1}{\rho_n}-\frac{1}{\rho}\right)\int_{r\ge \rho_n}\left(r\partial_ru_0\right)^2dr\right).\\
\end{eqnarray*}
Take $\rho_n$ sufficiently close to $\rho$ with $1-\frac{\rho_n}{\rho}<\frac{1}{4}$, then
\begin{equation}
\int_{|x|\ge \rho_n+t}|\nabla u^L|^2(x,t)+|\partial_t u^L|^2dx\ge \frac{1}{4}\int_{x\ge \rho_n}|\nabla u_0|^2+u_1^2dx, \,\,{\rm for}\,\,t\ge 0.
\end{equation}
Take $1>\delta>0$ small so that
\begin{equation}
\|a\|_{L^{5/4}_tL^{5/2}_x(|x|\leq \rho+1,\,0\leq t\leq \delta)}
\end{equation}
 is sufficiently small. Then we have
\begin{equation}
\int_{|x|\ge \rho_n+t}|\nabla(u^L-u)|^2+|\partial_t(u^L-u)|^2(t,x)dx\leq \frac{1}{100}\int_{x\ge \rho_n}|\nabla u_0|^2+u_1^2dx,\,\,{\rm for}\,\,0\leq t\leq \delta.
\end{equation}
Thus we obtain
\begin{equation}
\int_{|x|\ge \rho_n+t}|\nabla u|^2+(\partial_t u)^2(t,x)dx\ge \frac{1}{8}\int_{x\ge \rho_n}|\nabla u_0|^2+u_1^2dx>0,\,\,{\rm for}\,\,0\leq t\leq \delta.
\end{equation}
Thus $\rho(u(t),\partial_tu(t))\ge\rho_n+t$ for $0\leq t\leq \delta$. Passing $n$ to infinity, we obtain 
$$\rho(u(t),\partial_tu(t))\ge\rho(u_0,u_1)+t, \,\,{\rm for} \,\,0\leq t\leq \delta.$$
Then we can repeat this argument at $t=\delta$. The only danger is that the assumption there exists a sequence of radii approaching $\rho(\overrightarrow{u}(\delta))$
with the ``channel of energy" inequality valid for $t\ge \delta$ may fail. We observe that this can not happen, otherwise the support would 
expand backward and lead to $\rho(u(t),\partial_tu(t))>\rho+\delta$ for some $t<\delta$, a contradiction with finite speed propagation. The lemma is proved.\\

For later applications we also need the following fact for an elliptic problem. This result is certainly known to experts, however we are not
able to locate it in the literature, and thus include a proof in the Appendix for the convenience of readers.
\begin{lemma}\label{lm:elliptic}
Let $V\in Y$. For any $c\in R$, there exists a unique radial solution $u_c\in \dot{H}^1(B_r^c)$ for any $r>0$ to
\begin{equation}\label{eq:elliptic}
-\Delta u-V(x)u+u^5=0,\,\,{\rm in}\,\,R^3\backslash\{0\},
\end{equation}
with
\begin{equation}
\left|u(x)-\frac{c}{|x|}\right|=o(\frac{1}{|x|}),\,\,{\rm as}\,\,|x|\to\infty.
\end{equation}
If $u_c\in \dot{H}^1(R^3)$, then $u_c\in C^1(R^3)$ and
\begin{equation}
-\Delta u_c-Vu_c+u_c^5=0\,\,{\rm in}\,\,R^3.
\end{equation}
\end{lemma}

\smallskip
\noindent
Next we give a characterization of solutions to equation (\ref{eq:main}) that does not satisfy the channel of energy inequality outside some light cone. 
\begin{lemma}\label{lm:othersolutions}
Let $u\in C(R,\dot{H}^1)\cap L^5_tL^{10}_x((-T,T)\times R^3)$ for any $T>0$ be the solution to equation (\ref{eq:main}) with radial finite energy initial data
$\overrightarrow{u}(0)=(u_0,u_1)$. Suppose for some $R>0$ we have
\begin{equation}
\limsup_{|t|\to\infty}\int_{|x|\ge R+|t|}|\nabla u|^2+(\partial_tu)^2(x,t)dx=0.
\end{equation}
Then either $(u_0,u_1)$ is compactly supported, or for some $c\in R$, $(u_0,u_1)-(u_c,0)$ is compactly supported. 
\end{lemma}

\smallskip
\noindent
{\bf Remark.} As is observed in \cite{DKM} 
\begin{equation*}
\limsup_{t\to\infty}\int_{|x|\ge R+|t|}|\nabla u|^2+(\partial_tu)^2(x,t)dx>0
\end{equation*}
is equivalent to
\begin{equation*}
\inf_{t>0}\int_{|x|\ge R+|t|}|\nabla u|^2+(\partial_tu)^2(x,t)dx>0.
\end{equation*}
This follows from finite speed of propagation.\\

\smallskip
\noindent
{\bf Proof.} The proof is almost the same as in \cite{DKM}. The only difference is when $ru_0(r)\to l\neq 0$ (in their notation), we need to compare $\overrightarrow{u}$ with $(u_l,0)$ instead of 
a rescaled ground state $W$ (see also \cite{DKM2,KLS} for similar arguments). The possible slow decay of $V$ only slightly complicates the argument. However we note that the decay condition on $V$ we assume here seems to be sharp. For the convenience of readers we include the proof in Appendix B. Here we just emphasize the fact that we can find solutions to equation (\ref{eq:elliptic}) with asymptotic $\frac{c}{|x|}$ for any $c\in R$ is crucial for this argument. \\

The next theorem shows the only solutions that do not satisfy the ``channel of energy " inequality for all light cone are steady states.
\begin{theorem}\label{th:channelofenergy}
Suppose radial finite energy $(u_0,u_1)\not\equiv (u_c,0)$ for any steady state solution of equation (\ref{eq:mainequation}). Let $u\in C(R,\dot{H}^1)\cap L^5_tL^{10}_x((-T,T)\times R^3)$ for any $T\in(0,\infty)$ be the unique solution to equation (\ref{eq:mainequation}) with $\overrightarrow{u}(0)=(u_0,u_1)$. Then there exists $R>0$ such that
\begin{equation}\label{eq:channelofenergy}
\int_{|x|\ge R+|t|} |\nabla u|^2+(\partial_tu)^2(t,x)\,dx\ge \delta>0,
\end{equation}
for all $t\ge 0$ or all $t\leq 0$.
\end{theorem}

\smallskip
\noindent
{\bf Remarks.} In particular for the equation $\partial_{tt}u-\Delta u+u^5=0$ which has no nontrivial steady state, we have ``channel of energy" inequality (\ref{eq:channelofenergy}) for all nontrivial solutions.\\

\smallskip
\noindent
{\bf Proof.} By Lemma \ref{lm:othersolutions} we only need to consider the case when $(u_0,u_1)$ is compactly supported and the case $(u_0,u_1)-(u_c,0)$ is compactly supported for some steady state $(u_c,0)$. \\
Let us first consider the case $(u_0,u_1)$ is compactly supported. Then by Lemma \ref{lm:propagationofsupport} at $t_0$ (without loss of generality we take $t_0>0$) we have
\begin{equation}
\rho(u(t_0),\partial_tu(t_0))=\rho(u_0,u_1)+t_0.
\end{equation}
Take $t_0$ so large such that $\|V\|_{L^{5/4}_tL^{5/2}_x(\{(x,t)|\,|x|\ge\rho(u(t_0),\partial_tu(t_0))+t-t_0-1,\,t\ge t_0\})}$ is sufficiently small, and take $\rho\in(t_0,\rho(\overrightarrow{u}(t_0)))$ sufficiently close to $\rho(\overrightarrow{u}(t_0))$, then we have for $t\ge t_0$
\begin{equation}
\int_{|x|\ge t-t_0+\rho}|\nabla u^L|^2+(\partial_tu^L)^2(t,x)dx\ge \frac{1}{4}\int_{|x|\ge\rho}|\nabla u|^2(t_0,x)+(\partial_tu)^2(t_0,x)dx,
\end{equation}
and
\begin{equation}
\int_{|x|\ge t-t_0+\rho}|\nabla (u^L-u)|^2+(\partial_t(u^L-u))^2(t,x)dx\leq\frac{1}{8}\int_{|x|\ge\rho}|\nabla u|^2(t_0,x)+(\partial_tu)^2(t_0,x)dx,
\end{equation}
where $\overrightarrow{u}^L$ is the solution to free wave equation with $\overrightarrow{u}^L(t_0)=\overrightarrow{u}(t_0)$. Thus
\begin{equation}
\int_{|x|\ge t-t_0+\rho}|\nabla u|^2+(\partial_tu)^2(t,x)dx\ge \frac{1}{8}\int_{|x|\ge\rho}|\nabla u|^2(t_0,x)+(\partial_tu)^2(t_0,x)dx>0,\,\,{\rm for}\,\,t\ge t_0.
\end{equation}
Simple analysis involving the geometry of light cones and finite speed of propagation imply for $t\ge 0$ and some $\delta>0$
\begin{equation}
\int_{|x|\ge \rho+t-t_0}|\nabla u|^2+(\partial_tu)^2(t,x)dx\ge \delta.
\end{equation}
Now let us consider the case when $(u_0,u_1)-(u_c,0)$ is compactly supported. By Lemma \ref{lm:globalregularityofperturbedwave} and remarks below it, we have $u,\,u_c\in L^5_tL^{10}_x(|x|\ge |t|)$.  Let $h=u-u_c$, then $(h_0,h_1)=(u_0,u_1)-(u_c,0)$ is compactly supported and $h$ satisfies
\begin{equation}
\partial_{tt}h-\Delta h-Vh+5u_c^4h+10u_c^3h^2+10u_c^2h^3+5u_ch^4+h^5=0\,\,{\rm in}\,\, R\times R^3.
\end{equation}
Again by Lemma \ref{lm:propagationofsupport} without loss of generality we assume the support of $\overrightarrow{h}(t)$ expand in forward time. Denote 
\begin{equation}
a(x,t)=-V+5u_c^4+10u_c^3h+10u_c^2h^2+5u_ch^3+h^4.
\end{equation}
Note
\begin{equation}
\lim_{R\to\infty}\|a\|_{L^{5/4}_tL^{5/2}_x(|x|\ge R+|t|)}=0.
\end{equation}
Take $t_0$ large (so that the support $\rho(\overrightarrow{h}(t_0))=t_0+\rho(\overrightarrow{h}(0))$ is also large) and $\rho>t_0$ sufficiently close to $\rho(\overrightarrow{h}(t_0))$, such that by perturbation result in Lemma \ref{lm:perturbationresult} we have for $t\ge t_0$
\begin{equation}
\int_{|x|\ge t-t_0+\rho}|\nabla h^L|^2+(\partial_th^L)^2(t,x)dx\ge \frac{1}{4}\int_{|x|\ge\rho}|\nabla h|^2(t_0,x)+(\partial_th)^2(t_0,x)dx,
\end{equation}
and
\begin{equation}
\int_{|x|\ge t-t_0+\rho}|\nabla (h^L-h)|^2+(\partial_t(h^L-h))^2(t,x)dx\leq\frac{1}{8}\int_{|x|\ge\rho}|\nabla h|^2(t_0,x)+(\partial_th)^2(t_0,x)dx,
\end{equation}
where $\overrightarrow{h}^L$ is the solution to free wave equation with $\overrightarrow{h}^L(t_0)=\overrightarrow{h}(t_0)$. Thus
\begin{equation}
\int_{|x|\ge t-t_0+\rho}|\nabla h|^2+(\partial_th)^2(t,x)dx\ge \frac{1}{8}\int_{|x|\ge\rho}|\nabla h|^2(t_0,x)+(\partial_th)^2(t_0,x)dx>0,\,\,{\rm for}\,\,t\ge t_0.
\end{equation}
By the decay property of $(u_c,0)$ we obtain
\begin{equation}
\limsup_{t\to\infty}\int_{|x|\ge t-t_0+\rho}|\nabla u|^2+(\partial_tu)^2(t,x)dx\ge \frac{1}{8}\int_{|x|\ge\rho}|\nabla u|^2(t_0,x)+(\partial_tu)^2(t_0,x)dx>0.
\end{equation}
Thus we obtain by finite speed of propagation
\begin{equation}
\inf_{t>0}\int_{|x|\ge t-t_0+\rho}|\nabla u|^2+(\partial_tu)^2(t,x)dx>0.
\end{equation}

\smallskip
\noindent
{\bf Remarks.} Now let $\overrightarrow{U}^L$ be a nontrivial solution to free wave equation, and let sequence $\lambda_n>0$ satisfy $\lim\limits_{n\to\infty}(\lambda_n+\frac{1}{\lambda_n})=\infty$. Consider the nonlinear profile $U$ associated to $U^L,\,\lambda_n,\,t_n\equiv 0$ as in Lemma \ref{lm:aux}. Let $U_n(t,x)=\frac{1}{\lambda_n^{1/2}}U(\frac{t}{\lambda_n},\frac{x}{\lambda_n})$. By a simple change of variable we obtain from the remark below the above lemma for all $t\ge 0$ or all $t\leq 0$
\begin{equation}
\int_{|x|\ge |t|}|\nabla U_n|^2+(\partial_tU_n)^2(t,x)dx=\int_{|x|\ge |t|}|\nabla U|^2+(\partial_tU)^2(t,x)dx>\delta>0.
\end{equation} 

\medskip
\noindent
We shall also need the following ``channel of energy" inequality for profiles with $\lim\limits_{n\to\infty}\frac{t_n}{\lambda_n}\in\{\pm\infty\}$.
\begin{lemma}\label{lm:channelforprofile}
Let $U^L_{jn}(x,t)=\frac{1}{\lambda_{jn}^{1/2}}U^L_j(\frac{t-t_{jn}}{\lambda_{jn}},\frac{x}{\lambda_{jn}})$ be with 
\begin{equation}
\lim_{n\to\infty}\frac{t_{jn}}{\lambda_{jn}}=-\infty \,\,({\rm or}\,\,+\infty).
\end{equation}
Let $U_j$ be the nonlinear profile associated to $U^L_j,\,\lambda_{jn},\,t_{jn}$ as in Lemma \ref{lm:aux}. Then for $n$ sufficienly large and all $t>0$ (or all $t<0$ respectively), we have
\begin{equation}\label{lm:channelforprofile} 
\int_{|x|\ge t}|\nabla U_{jn}|^2+(\partial_t U_{jn})^2(t,x)dx\ge \epsilon_0>0,
\end{equation}
where $U_{jn}(x,t)=\frac{1}{\lambda_{jn}^{1/2}}U_j(\frac{t-t_{jn}}{\lambda_{jn}},\frac{x}{\lambda_{jn}})$.
\end{lemma}

\smallskip
\noindent
{\bf Proof.} By the definition of $\overrightarrow{U}^j$,
\begin{equation}
\lim_{t\to\infty}\|\overrightarrow{U}_j(t)-\overrightarrow{U}_j^L(t)\|_{\dot{H}^1\times L^2}=0.
\end{equation}
By Lemma \ref{lm:channelwithdirection}, there exists $t_0>0$ such that for $t\ge t_0$
\begin{eqnarray*}
&&\frac{1}{|S^2|}\int_{|x|\ge t-t_0}|\nabla U^L_j|^2+(\partial_tU^L_j)^2(t,x)dx\\
&&\quad\ge\int_{r\ge t-t_0}\left(\partial_r(rU_j^L)\right)^2+(\partial_t rU_j^L)^2(t,r)dr\\
&&\quad\ge \frac{1}{2}\int_{r>0}\left(\partial_r(rU_j^L)\right)^2+(\partial_t rU_j^L)^2(t_0,r)dr\ge \frac{\epsilon_0}{|S^2|}>0.
\end{eqnarray*}
Thus for some $t_1$ sufficiently large we have for $t\ge t_1$
\begin{equation}
\int_{|x|\ge t-t_0}|\nabla U_j|^2(t,r)+(\partial_t U_j)^2(t,x)dx\ge \frac{\epsilon_0}{2}.
\end{equation}
Thus by a simple change of variables, we see for $t-t_{jn}\ge\lambda_{jn}t_1$,
\begin{equation}
\int_{|x|\ge t-(t_{jn}+\lambda_{jn}t_0)} |\nabla U_{jn}|^2(r,t)+(\partial_t U_{jn})^2(t,x)dx \ge\frac{\epsilon_0}{2}.
\end{equation}
By the assumptions on $\lambda_{jn},\,t_{jn}$, for $n$ sufficiently large $t_{jn}+\lambda_{jn}t_1\leq 0$ and $t_{jn}+\lambda_{jn}t_0\leq 0$. The lemma is proved.\\

\smallskip
\noindent
In summary, the above lemmas show that in the profile decompositions (Lemma \ref{lm:profilewithpotential}) only the profile given by steady state solution does not satisfy the following ``channel of energy" inequality:
\begin{equation}\label{eq:nonlinearchannel}
\int_{|x|\ge |t|}|\nabla U_{jn}|^2+(\partial_t U_{jn})^2(t,x)dx\ge \epsilon_0>0,
\end{equation}
for $n$ sufficiently large and all $t\ge 0$ or all $t\leq 0$.\\

\begin{section}{Proof of main result}
\begin{theorem}
Suppose $(u_0,u_1)\in \dot{H}^1\times L^2$ is radial, and $V\in Y$.  Denote the set of radial steady states of equation (\ref{eq:main}) as $\Sigma$. Let $u\in C([0,\infty),\dot{H}^1)\cap L^5_tL^{10}_x([0,T)\times R^3)$ for any $T<\infty$ be the unique solution to equation (\ref{eq:main}) with initial data $(u_0,u_1)$. Then there
exists a solution $\overrightarrow{u}^L$ to the free wave equation, such that 
\begin{equation}
\lim_{t\to\infty}\inf_{(u_c,0)\in\Sigma}\|\overrightarrow{u}(t)-\overrightarrow{u}^L(t)-(u_c,0)\|_{\dot{H}^1\times L^2}=0.
\end{equation}
Moreover, there exists a dense open subset $\Omega\subseteq Y$, such that if $V\in Y$ then there are only finitely many radial steady states for equation (\ref{eq:main}). In this case there exists some steady state $(u_c,0)$ so that we have
\begin{equation}
\lim_{t\to\infty}\|\overrightarrow{u}(t)-\overrightarrow{u}^L(t)-(u_c,0)\|_{\dot{H}^1\times L^2}=0.
\end{equation}
\end{theorem}

\smallskip
\noindent
{\bf Proof.} Take a sequence of time $t_n\to\infty$, and let $(u_{0n},u_{1n})=(u(t_n),\partial_tu(t_n))$. Passing to a subsequence if necessary, we can assume $(u_{0n},u_{1n})$ has the following profile decomposition
\begin{equation}
(u_{0n},u_{1n})=\overrightarrow{u}^L(t_n)+\overrightarrow{U}^L_1(0)+\sum_{j=2}^J(\frac{1}{\lambda_{jn}^{1/2}}U^L_j(-\frac{t_{jn}}{\lambda_{jn}},\frac{x}{\lambda_{jn}}),\frac{1}{\lambda_{jn}^{3/2}}\partial_tU^L_j(-\frac{t_{jn}}{\lambda_{jn}},\frac{x}{\lambda_{jn}}))+\overrightarrow{w}_{Jn}(0),
\end{equation}
where $\overrightarrow{u}^L$ is the free wave for $\overrightarrow{u}$ given by Theorem \ref{th:freewave}, and the parameters $\lambda_{jn},\,t_{jn}$ satisfy the usual orthogonality property. Our potential $V$ does not decay in time in the whole space. Thus at first glance we can not apply Lemma \ref{lm:profilewithpotential}. Note however in the exterior of light cone $S=\{(x,t)|\,|x|\ge |t|\}$, $V$ satisfies the bound $\|V\|_{L^{5/4}_tL^{5/2}_x\cap L^1_tL^3_x(S)}<\infty$ required in Lemma \ref{lm:profilewithpotential}. Therefore by Lemma \ref{lm:profilewithpotential} and finite speed of propagation, the solution $\overrightarrow{u}_n(t)$ to equation (\ref{eq:main}) with $\overrightarrow{u}_n(0)=(u_{0n},u_{1n})$ has the following decomposition for $|x|\ge |t|$
\begin{equation}
\overrightarrow{u}_n(t)=\overrightarrow{u}^L(t+t_n)+\overrightarrow{U_1}+\sum_{j=2}^J \overrightarrow{U}_{jn}+\overrightarrow{w}_{Jn}+\overrightarrow{r}_{Jn},\,\,
\end{equation}
with 
\begin{equation}
\lim_{J\to\infty}\limsup_{n\to\infty}\sup_{t\in R}\|r_{Jn}(t)\|_{\dot{H}^1\times L^2}=0,
\end{equation}
where $\overrightarrow{U}_{jn}(x,t)=(\frac{1}{\lambda_{jn}^{1/2}}U_j(\frac{t-t_{jn}}{\lambda_{jn}},\frac{x}{\lambda_{jn}}),\frac{1}{\lambda_{jn}^{3/2}}\partial_tU_j(\frac{t-t_{jn}}{\lambda_{jn}},\frac{x}{\lambda_{jn}}))$. Moreover, denoting 
\begin{equation}
\overrightarrow{U}_{0n}(t)=\overrightarrow{u}^L(t+t_n),\,\,{\rm and}\,\,\overrightarrow{U}_{1n}=\overrightarrow{U}_1,
\end{equation} 
then for $\rho_n>\sigma_n>0$ and $\theta_n\in R$ we have the following orthogonality property for $0\leq j\not=j'$
\begin{eqnarray}
&&\lim_{n\to\infty}\int_{\sigma_n<|x|<\rho_n}\nabla U_{jn}\nabla U_{j'n}+\partial_tU_{jn}\partial_tU_{j'n}(\theta_n,x)dx=0;\\
&&\lim_{n\to\infty}\int_{\sigma_n<|x|<\rho_n}\nabla U_{jn}\nabla w_{Jn}+\partial_tU_{jn}\partial_tw_{Jn}(\theta_n,x)dx=0.\\
\end{eqnarray}
We note that we can leave the profile $\overrightarrow{u}^L$ linear, since the nonlinear profile is close to this linear profile as $n\to\infty$ in $\{|x|\ge|t|\}$, and the error term is absorbed in $r_{Jn}$. It should also be emphasized that $U_1$ solves equation (\ref{eq:main}) and $U_j$ solves nonlinear wave equation without potential for $j\ge 2$, as in Lemma \ref{lm:profilewithpotential}. 
Suppose for some $j\ge 1$ and all sufficiently large $n$ we have for all $t\ge 0$
\begin{equation}
\int_{|x|\ge t}|\nabla U_{jn}|^2+(\partial_tU_{jn})^2(t,x)dx>\delta>0.
\end{equation}
Then by the orthogonality property in Lemma \ref{lm:profilewithpotential} and the error estimate on $r_{Jn}$, we obtain for $n$ sufficiently large and $t\ge 0$
\begin{equation}
\|\overrightarrow{u}_n(t)-u^L(t+t_n)\|^2_{\dot{H}^1\times L^2(|x|\ge |t|)}\ge \frac{1}{2}\int_{|x|\ge |t|}|\nabla U_{jn}|^2+(\partial_tU_{jn})^2(x,t)dx.
\end{equation}
Since $\overrightarrow{u}_n(t,x)=\overrightarrow{u}(t+t_n,x)$ for $|x|\ge |t|$ and by Theorem \ref{th:freewave} for any fixed $n$
\begin{equation}
\lim_{t\to\infty}\|\overrightarrow{u}(t+t_n)-\overrightarrow{u}^L(t+t_n)\|_{\dot{H}^1\times L^2(|x|\ge t)}=0.
\end{equation}
We obtain for $n$ sufficiently large
\begin{equation}
\limsup_{t\to\infty}\int_{|x|\ge t}|\nabla U_{jn}|^2+(\partial_tU_{jn})^2(t,x)dx=0,\,\,{\rm a \,\,contradiction}.
\end{equation}

Now suppose for some $j\ge 1$ and sufficiently large $n$ we have 
\begin{equation}
\inf_{t<0}\int_{|x|\ge |t|}|\nabla U_{jn}|^2+(\partial_tU_{jn})^2(t,x)dx>0.
\end{equation}
Applying the orthogonality property at $\theta_n=-t_n$ we would get for $n$ sufficiently large
\begin{equation}
\int_{|x|\ge t_n}|\nabla u_0|^2+u_1^2(x)dx\ge \frac{1}{2}\inf_{t<0}\int_{|x|\ge |t|}|\nabla U_{jn}|^2+(\partial_tU_{jn})^2(t,x)dx>0.
\end{equation}
A contradiction with $(u_0,u_1)\in \dot{H}^1\times L^2$. Thus all $U_{jn}$ have no ``channel of energy" property for all $1\leq j\leq J$. By results in Section 4, we see that $U_j\equiv 0$ for all $j\ge 2$ and $U_1=u_c$ for some steady state. Thus there are at most $2$ profiles in the decomposition. We thus obtain for $J_0= 2$
\begin{equation}
\limsup_{n\to\infty}\|w_{J_0n}\|_{L^5_tL^{10}_x(R\times R^3)}=0.
\end{equation}
Recall that $\overrightarrow{w}_{J_0n}$ solves the free wave equation. We claim
\begin{equation}
\limsup_{n\to\infty}\|\overrightarrow{w}_{J_0n}(t)\|_{\dot{H}^1\times L^2}=0.
\end{equation}
Otherwise by passing to a subsequence if necessary we would have by channel of energy inequality for linear wave equation in Lemma \ref{lm:channelofenergyforlinearwave}, that for sufficiently large $n$ and all $t\ge 0$ or all $t\leq 0$,
\begin{equation}
\int_{|x|\ge |t|}|\nabla w_{J_0n}|^2+(\partial_tw_{J_0n})^2(t,x)dx\ge \delta>0.
\end{equation}
This leads to contradiction by similar arguments as in the last paragraph using orthogonality property. Thus we obtain
\begin{equation}
(u(t_n),\partial_tu(t_n))=\overrightarrow{u}^L(t_n)+(u_c,0)+\overrightarrow{w}_{J_0n}(0)
\end{equation}
and $\lim\limits_{n\to\infty}\|\overrightarrow{w}_{J_0n}(t)\|_{\dot{H}^1\times L^2}=0$. Thus 
\begin{equation}\label{eq:convergencetosets}
\lim_{n\to\infty}\inf_{(u_c,0)\in\Sigma}\|\overrightarrow{u}(t_n)-\overrightarrow{u}^L(t_n)-(u_c,0)\|_{\dot{H}^1\times L^2}=0.
\end{equation}
Since this is true for some subsequence of an arbitrary sequence of time going to infinity the first part of theorem is proved. By Theorem \ref{th:numberofsteadystates} below, there exists a dense open subset $\Omega\subset Y$, such 
that if $V\in \Omega$ then there are only finitely many radial steady states for equation (\ref{eq:main}). In this case the second part follows immediately from the first part of the theorem.

\end{section}

\begin{section}{Appendix A. The elliptic problem}
In this Appendix we first give a proof of Lemma \ref{lm:elliptic}.\\
{\bf Main result}\\
Let $V\in Y$. For any $c\in R$, there exists a unique radial solution
$u\in \dot{H}^1(B_r^c)$ for any $r>0$ to
\begin{equation}
-\Delta u-V(x)u+u^5=0,\,\,{\rm in}\,\,R^3\backslash\{0\},
\end{equation}
with the asymptotics
\begin{equation}
\left|u(x)-\frac{c}{|x|}\right|=o(\frac{1}{|x|}), \,\,{\rm as}\,\,|x|\to\infty.
\end{equation}
If $u_c\in \dot{H}^1(R^3)$, then $u_c\in C^1(R^3)$ and
\begin{equation}
-\Delta u_c-V(x)u_c+u_c^5=0,\,\,{\rm in}\,\, R^3.
\end{equation}

\bigskip
\noindent
We first prove the following lemmas, from which our main result follows easily.
\begin{lemma}\label{lm:boundaryvalueproblem}
Let $V$ be as above. Suppose $R$ is sufficiently large such that
\begin{equation}\label{eq:smallnessofpotential}
\int_{B_R^c}V(x)u^2(x)dx\leq\frac{1}{2}\int_{B_R^c}|\nabla u|^2dx
\end{equation}
holds for all $u\in\dot{H}^1(B_R^c)$. Consider the following boundary value problem
\begin{equation}\label{eq:boundaryvalueproblem}
-\Delta u-V(x)u+u^5=0,\,\,{\rm in}\,\,B_R^c,
\end{equation}
with boundary value $u|_{\partial B_R}=\lambda$. Then for any $\lambda\in R$ there exists a unique solution $u\in\dot{H}^1(B_R^c)$.
\end{lemma}

\smallskip
\noindent
{\bf Remarks.} Although $V$ may change sign, the condition (\ref{eq:smallnessofpotential}) which is essentially a ``smallness" condition on the positive part $V^{+}$ of potential, implies that we will have the usual ``comparison principle"
valid (see below for details). This is of course well known and we will make use of this observation below.

\smallskip
\noindent
{\bf Proof.} We only need to consider the case $\lambda\ge 0$ by the invariance of the equation under the transformation $u\to -u$. Let us first prove existence of solution.
Fix $\lambda\ge 0$, consider the minimization problem
\begin{equation}
E(\lambda):=\inf_{u\in\dot{H}^1(B_R^c),\,u|_{\partial B_R}=\lambda}\,\,\int_{B_R^c}\frac{|\nabla u|^2}{2}-\frac{Vu^2}{2}+\frac{u^6}{6}dx.
\end{equation}
Note that our nonlinearity is defocusing and in this case, it's easy to show there is a minimizer of the above minimization problem (and the minimizer enjoys strong a priori estimates). It's then clear the minimizer satisfies the Euler-Lagrange equation of the functional which is exactly equation (\ref{eq:boundaryvalueproblem}) with the boundary value $u|_{\partial B_R^c}=\lambda$. Thus the existence is established. Next we observe that $u$ is nonnegative. This follows from multiplying the equation (\ref{eq:boundaryvalueproblem}) with $u$ and integrate over the region $\Omega:=\{x\in B_R^c:u(x)\leq 0\}$. With the observation that there are no boundary terms in the integration by parts formula, this leads to
\begin{equation}
0=\int_{\Omega}|\nabla u|^2-Vu^2+u^6dx\ge \int_{\Omega}\frac{|\nabla u|^2}{2}+u^6dx.
\end{equation}
Thus $u$ can never be negative. Then by strong maximum principle we have $u>0$ if $\lambda>0$. Now take solutions $u_1,\,u_2$ with boundary values $\lambda_1<\lambda_2$ respectively. Let $w=u_1-u_2$, then $w$ satisfies
\begin{equation}
-\Delta w-Vw+b(x)w=0,
\end{equation}
where $b(x)$ is nonnegative. Multiply the above equation with $w$ and integration by parts in the region $\Omega=\{w>0\}$, again there is no boundary term. Thus we obtain 
\begin{equation}
0=\int_{B_R^c}|\nabla w|^2-Vw^2+b(x)w^2\,dx\ge \frac{1}{2}\int_{B_R^c}|\nabla w|^2+b(x)w^2\,dx.
\end{equation}
Thus $\Omega=\emptyset$. By strong maximum principle we see $w<0$, that is $u_1< u_2$. Uniqueness follows immediately. The lemma is proved.\\

\smallskip
\noindent
The next lemma gives more precise asymptotics for the solutions obtained in the above lemma.
\begin{lemma}\label{lm:preciseasymptotics}
Let $V,\,R$ be as in the above lemma. Let $\lambda\ge 0$ and $u(\cdot,\lambda)\in \dot{H}^1(B_R^c)$ be the unique solution to equation (\ref{eq:boundaryvalueproblem}) with the boundary value $u|_{\partial B_R}=\lambda$. Then there exists $c\ge 0$ with
\begin{equation}\label{eq:asymptotics}
\left|u(x,\lambda)-\frac{c}{|x|}\right|=o(\frac{1}{|x|}), \,\,{\rm as}\,\, |x|\to\infty,
\end{equation} 
and we have the bound
\begin{equation}
|u(x,\lambda)|\leq\frac{C(\lambda,\|V\|_Y,R)}{|x|}\,\,{\rm for}\,\,|x|\ge R.
\end{equation}
Moreover if we denote $c=c(\lambda)$, then $c$ is continuous, strictly increasing with $\lambda$, and 
\begin{equation}
\lim\limits_{\lambda\to\infty}c(\lambda)=\infty.
\end{equation}
\end{lemma}

\smallskip
\noindent
{\rm Proof.} We first establish the asymptotics (\ref{eq:asymptotics}). Take the function $\phi$ with
\begin{equation}
-\Delta \phi=\frac{M}{(1+|x|^2)^{\frac{\beta}{4}+1}},\,\,{\rm in}\,\,R^3.
\end{equation}
Note that $\frac{\beta}{2}+2>3$. Then $\phi(x)$ is positive and comparable with $\frac{1}{|x|}$ as $x\to\infty$. Thus
\begin{equation}
-\Delta \phi(x)-V(x)\phi(x)+\phi^5(x)\ge \frac{M}{2|x|^{\frac{\beta}{2}+2}}-\frac{C}{|x|^{\beta}}\frac{CM}{|x|}+\phi^5(x)\ge 0,\,\,{\rm for}\,\,|x|\ge R_1>0,
\end{equation}
where $R_1=R_1(\lambda,\|V\|_Y,R)$.
Thus by comparison principle when $M=M(\lambda,\|V\|_Y,R)$ is sufficiently large we have
\begin{equation}
u(x)\leq\phi(x)\leq \frac{CM}{|x|},\,\,{\rm for}\,\,|x|\ge R_1.
\end{equation}
Once we have this decay estimate the asymptotics (\ref{eq:asymptotics}) follow from the property of Green's function for the Laplace in the exterior of a ball. Moreover by uniqueness and a priori bound of solutions, we can obtain the continuity of $u(x,\lambda)$ in $x\in K\Subset \overline{B_R^c}$ and $\lambda$. This, combined with the decay estimate easily implies the continuity of $c(\lambda)$ in $\lambda$.
Thus $c=c(\lambda)$ is well defined,  continuous, and by comparison principle is nondecreasing in $\lambda$. We still need to show it's strictly increasing in $\lambda$, and 
$\lim\limits_{\lambda\to\infty}c(\lambda)=\infty$.

\end{section}
We first prove $c(\lambda)$ is strictly increasing. Take $u_1$ and $u_2$ solutions to equation (\ref{eq:boundaryvalueproblem}) with boundary data $\lambda_1$ and $\lambda_2$ respectively. Assume $\lambda_1<\lambda_2$, let 
$w=u_2-u_1$, then $w>0$ and satisfies
\begin{equation}
-\Delta w-V(x)w+5u_1^4w+10u_1^{3}w^2+10u_1^{2}w^3+5u_1w^4+w^5=0\,\,{\rm in}\,\,B_R^c,
\end{equation}
with boundary value $w|_{\partial B_R}=\lambda_2-\lambda_1$. Denote $\beta_1=\min\{\beta,4\}$. We have the following estimates on the coefficient $b=-V+5u_1^4+10u_1^3w+10u_1^2w^2+5u_1w^3+w^4$:
\begin{equation}
|b(x)|\leq \frac{C(\lambda_2,\|V\|_Y,R)}{|x|^{\beta_1}}\,\,{\rm for}\,\,x\in B_R^c.
\end{equation}
Take $\alpha$ large, and let $v_{\alpha}(x)=\alpha^{\frac{1}{2}} w(\alpha x)$, $b_{\alpha}(x)=\alpha^2b(\alpha x)$. Then $v_{\alpha}$ satisfies
\begin{equation}
-\Delta v_{\alpha}+ b_{\alpha}(x)v_{\alpha}=0\,\,{\rm for}\,\, x\in B_R^c,
\end{equation}
and $v_{\alpha}|_{\partial B_R}>0$.
Note that 
\begin{equation}
|b_{\alpha}(x)|\leq\frac{C(\lambda_2,\|V\|_Y,R)}{\alpha^{\beta_1-2}|x|^{\beta_1}}\,\,{\rm for}\,\, |x|\ge R.
\end{equation} 
Thus if we take $\alpha$ sufficiently large, then $b_{\alpha}$ will be very small, so that we can use standard perturbation argument to solve $v_{\alpha}$ with asymptotics $\frac{c}{|x|}$ for some positive $c$ (in fact $v_{\alpha}$ is close to $\frac{v_{\alpha}(R)R}{|x|}$). Thus by the relation between $v_{\alpha}$ and $w$, we see $w(x)\sim \frac{\epsilon}{|x|}$ as $|x|\to\infty$ for some $\epsilon>0$. Thus $c(\lambda)$ is strictly increasing.\\
Lastly we prove that $\lim\limits_{\lambda\to\infty}c(\lambda)=\infty$. Note that for any fixed large $\alpha>0$, $u(,\lambda)|_{R\alpha\leq |x|\leq 2R\alpha}$ can not stay bounded as $\lambda\to\infty$. Otherwise by elliptic estimates we would have the gradient of $u(,\lambda)$ stays bounded in $\{\frac{4R\alpha}{3}\leq|x|\leq \frac{5R\alpha}{3}\}$. Since $u(,\lambda)$ is radial and thus satisfies an ODE in $r$. This would imply $u(,\lambda)|_{\partial B_R}=\lambda$ stays bounded as $\lambda\to\infty$, a contradiction. Assume for $\lambda$ sufficiently large at $|x|=R\alpha_{\ast}\in(R\alpha,2R\alpha)$ we have $u(x,\lambda)\ge 1$. Now define $v_{\alpha_{\ast}}(x)=\alpha_{\ast}^{\frac{1}{2}}u(\alpha_{\ast} x,\lambda)$ and $V_{\alpha_{\ast}}(x)=\alpha_{\ast}^2V(\alpha_{\ast} x)$, we have
\begin{equation}
-\Delta v_{\alpha_{\ast}}-V_{\alpha_{\ast}}v_{\alpha_{\ast}}+v_{\alpha_{\ast}}^5=0,\,\,x\in B_R^c.
\end{equation} 
with $v_{\alpha_{\ast}}|_{\partial B_R}\ge 1$. Since $V_{\alpha_{\ast}}$ is small, we have $v_{\alpha_{\ast} }(x)\ge\frac{\epsilon}{|x|}$ for $|x|\ge R$ by the following Lemma \ref{lm:smallresult}. Thus by a change of variable we get $u(x,\lambda)\ge \frac{\alpha_{\ast}^{1/2}\epsilon}{|x|}$ for $|x|\ge R\alpha_{\ast}$. This is true for sufficiently large $\lambda$ depending on the value of $\alpha$. Thus 
\begin{equation}
\liminf_{\lambda\to\infty}c(\lambda)\ge \alpha^{1/2}_{\ast}\epsilon\ge\alpha^{1/2}\epsilon.
\end{equation}
Since $\alpha$ is arbitrary. The proof is complete.\\
 
\smallskip
\noindent
We now finish the proof of the above lemma by stating the following perturbation result (without proof).
\begin{lemma}\label{lm:smallresult}
Let $R,\,\beta_1$ be defined as above. There exists $\epsilon,\,\delta>0$, such that if 
\begin{equation}
\sup_{x\in B_R^c}|x|^{\beta_1}|b(x)|<\delta,
\end{equation}
then any solution $u$ to
\begin{equation}
-\Delta u+b(x)u+u^5=0,\,\,x\in B_R^c
\end{equation}
with $u|_{\partial B_R}\ge 1$ satisfies
\begin{equation}
|u(x)|\ge \frac{\epsilon}{|x|},\,\,{\rm for}\,\,x\in B_R^c.
\end{equation}
\end{lemma}

\medskip
\noindent
{\bf Proof of Main result.} For any $c\in R$, take $\lambda$ with $c(\lambda)=c$ and let $u_c$ be the solution to equation (\ref{eq:boundaryvalueproblem}). Since $u_c$ is radial, we can formulate 
the equation as an ODE in $r$ for $r\ge R$. We can then solve this ODE from $R$ to $0$ and in this way extend $u_c$ as a solution in $R^3\backslash\{0\}$. It is easy to check this $u_c$ is the only solution satisfying conditions in Main Result. If $u\in \dot{H}^1(R^3)$ the other conclusions are well known from classical elliptic theory.

\bigskip
\noindent
The next result shows that for ``generic" potential $V$ there are only finitely many steady states. We have removed radial symmetry in our considerations below to have more general results. In the symmetric case the conclusions and proofs are exactly the same.
\begin{theorem}\label{th:numberofsteadystates}
Fix $\beta>2$. Define
\begin{eqnarray}
&&X:=\{\phi\in\dot{H}^1\cap L^6\cap\dot{H}^2(R^3)\},\\
&&{\rm with\,\,the\,\,natural\,\,norm\,\,}\|\phi\|_Y:=\|\phi\|_{\dot{H}^1}+\|\phi\|_{\dot{H}^2};\\
&&Y:=\{\phi\in C(R^3):\,\,\sup_{x\in R^3}((1+|x|)^{\beta}|\phi(x)|)<\infty\},\,\,\\
&&{\rm with\,\,the\,\,natural\,\,norm}\,\,\|\phi\|_Y:=\sup_{x\in R^3}((1+|x|)^{\beta}|\phi(x)|);\\
&&Z:=\{\phi\in\dot{H}^{-1}\cap L^2(R^3)\},\\
&&{\rm with\,\,the\,\,natural\,\,norm\,\,}\|\phi\|_Z:=\|\phi\|_{\dot{H}^{-1}}+\|\phi\|_{L^2}.
\end{eqnarray}
Define 
\begin{equation}
F:X\times Y\to Z,\,\,{\rm with}\,\,F(u,V)=-\Delta u-Vu+u^5.
\end{equation}
Then there exists a dense open set $\Omega\subset Y$ such that for any $V\in \Omega$ the elliptic equation
\begin{equation}
F(u,V)=0
\end{equation}
has only finitely many solutions $u\in X$.
\end{theorem}

\smallskip
\noindent
{\bf Proof.} It's not hard to check that for $u\in X$ the linear operator
\begin{equation}
\mathcal{L}(u,V):X\to Z\,\,{\rm with}\,\,\mathcal{L}(u,V)\phi=-\Delta\phi-V\phi+5u^4\phi
\end{equation}
is a compact perturbation of the operator $-\Delta:X\to Z$. It is also clear that $-\Delta$ is invertible as an operator from $X\to Z$. 
Take a dense open subset $\Omega_1$ of $Y$ such that for any $V\in \Omega_1$ $-\Delta-V$ has no zero resonance or eigenvalue. It is not hard to prove that any $\phi\in\,{\rm ker}\,(-\Delta-V)$ is of order $|\phi(x)|\sim O(\frac{1}{|x|})$ as $|x|\to\infty$, using arguments similar to those in the proof of Lemma \ref{lm:preciseasymptotics}. Thus
when $V\in \Omega_1$, $-\Delta-V$ is invertible. We now apply the ``transversality theorem" from page 295 in \cite{ST} to $F:X\times \Omega_1\to Z$. The
main condition to check is the ``nondegeneracy" condition. More precisely we need to verify that the full linearization of $F$ at a solution $(u,V)$ to $F(u,V)=0$ is
surjective. We have
\begin{equation}
dF(u,V)\,(\phi,a(x))=-\Delta \phi-V\phi+5u^4\phi-a(x)u,
\end{equation}
where $\phi\in X,\,a\in Y$. We need to show for any $f\in Z$ there is $(\phi,a)\in X\times Y$ such that 
\begin{equation}\label{eq:existence}
dF(u,V)\,(\phi,a(x))=f. 
\end{equation}
We now show we can always solve equation (\ref{eq:existence}). We first consider the case $u$ is not identically zero. It suffices to consider the case when the operator $-\Delta-V+5u^4$ has nontrivial kernel. Denote
\begin{equation}
{\rm ker}\,(-\Delta-V+5u^4)=\,{\rm span}\,\{\psi_j,\,j=1,2,\dots,k\},
\end{equation}
where $\psi_j\in X$ are linearly independent. Let us define linear functionals $\Lambda_j$ on $Z$ for $j=1,2\dots,k$ as follows
\begin{equation}
\Lambda_j(g):=\int_{R^3}\psi_j(x)g(x)dx,
\end{equation}
firstly for $g\in C_c(R^3)$. Note we have the obvious bound $|\Lambda_j(g)|\leq C\|\psi_j\|_{\dot{H}^1}\|g\|_{\dot{H}^{-1}}$, thus we can extend this functional boundedly to $Z$. It is not hard to verify that 
\begin{equation}
{\rm Range}\,(-\Delta-V+5u^4)=\{g\in Z:\,\Lambda_j(g)=0,\,\,{\rm for\,\,all}\,\,j\leq k\}.
\end{equation}
Thus to prove the solvability of equation (\ref{eq:existence}), we only need to show that we can find $a\in Y$ such that 
\begin{equation}
\int_{R^3}a(x)u(x)\psi_j(x)dx=-\Lambda_j(f),\,\,{\rm for\,\,all\,\,}j\in\{1,2,\dots,k\}.
\end{equation}
Since $\{u(x)\psi_j(x),\,1\leq j\leq k\}$ is linearly independent by well known unique continuation results, the existence of $a$ solving the above equations is obvious. Thus equation (\ref{eq:existence}) has a solution $(\phi,a)$.
If $u\equiv 0$, then by our assumption on $\Omega_1$, $dF(0,V)=-\Delta -V$ is surjective. Thus in summary, the nondegeneracy condition in \cite{ST} is verified, and we can conclude that there exists a dense subset $\Omega$ of $\Omega_1$ such that for $V\in\Omega$ the solutions of $F(u,V)=0$ are simple, meaning that $F_u(u,V)\phi=-\Delta\phi-V\phi+5u^4\phi$ has trivial kernel. From this and the strong a priori estimates for solutions, we easily deduce that for $V\in\Omega$, there are only finitely many steady states.\\

\smallskip
\noindent
As is already mentioned in the introduction, in general we expect many ``excited states" for our equation (\ref{eq:main}). Here we would like to illustrate the mechanism for the creation of ``excited states" in a very simple case (in particular we assume $\beta$ to be large to avoid technical issues from now on). If we use the notation as in the above theorem and consider the solutions to $F(u,\alpha V)=0$ with a parameter $\alpha>0$. Then we increase the parameter $\alpha$ from values. Roughly speaking one pair of ``excited states" ($\phi$ and $-\phi$) will be created from zero solution whenever a negative eigenvalue of $-\Delta -\alpha V$ appears through standard bifurcation of pitchfork type. The picture is quite dynamical with respect to $\alpha$. 
\begin{theorem}\label{th:excitedstates}
Suppose $V\in Y$ is nonnegative and not identically zero. Suppose further the (Birman-Schwinger) compact operator $\sqrt{V}(-\Delta)^{-1}\sqrt{V}$ from $L^2$ to itself has only simple eigenvalues. 
Denote $\lambda_1,\,\lambda_2$ as the first and second largest eigenvalues of $\sqrt{V}(-\Delta)^{-1}\sqrt{V}$ respectively. Then for small $\epsilon,\,\delta>0$ and $\alpha\in(\frac{1}{\lambda_2},\frac{1}{\lambda_2}+\delta)$, there are exactly two nontrivial solutions $\phi(\alpha)$ and $-\phi(\alpha)\in B_{\epsilon}(0)\subset X$ to $F(u,\alpha V)=0$. Moreover $\phi(\alpha)$ changes sign.
\end{theorem}

\smallskip
\noindent
{\bf Proof.} Observing that at $\alpha=\frac{1}{\lambda_2}$ the operator $-\Delta-\alpha V$ has a simple kernel. Then the existence and uniqueness follows immediately from standard bifurcation theorems. The claim that these solutions should change sign follows from the fact that they are close to a scalar multiple of functions in the kernel which changes sign since the smallest eigenvalue of $-\Delta-\alpha V$ is negative. These steady states are clearly ``excited states".\\

\smallskip
\noindent
{\bf Remarks.} In general for radial $V$ if $\sqrt{V}(-\Delta)^{-1}\sqrt{V}$ has $n$ simple eigenvalues $> 1$ with radial eigenfunctions. Then by using the methods in \cite{YanYan} one should be able to show there are at least $2n+1$ steady states to equation (\ref{eq:main}). Such methods are based on global bifurcation theory of Rabinowitz and a ``counting of number of zeroes" argument to distinguish different bifurcation curves. We will however not pursue this direction here.\\

\smallskip
\noindent
Suppose $V$ is generic in the sense of Theorem \ref{th:numberofsteadystates} and $(\psi(x),0)$ is an ``excited state". Thus the linearized operator around $\psi$, $-\Delta-V+5\psi^4$ has no zero resonance or zero eigenvalues (so that we have the necessary dispersive estimates for the continuous part of the linearized operator). In such situations it's by now standard how to construct ``center stable" manifolds near $(\psi,0)$. In particular there are solutions $\overrightarrow{u}(t)$ to equation (\ref{eq:main}) which scatters to the excited state $(\psi,0)$. We refer readers to \cite{SchNak} for details. These discussions demonstrate in general the global dynamical picture for equation (\ref{eq:main}) can be very complicated.

\end{section}

\begin{section}{Appendix B.}
Recall the main equation
\begin{equation}
\label{CP} \left\{
\begin{gathered}
\partial_t^2 u -\Delta u-V(x)u+u^5=0,\quad (t,x)\in \R\times \R^3\\
u_{\restriction t=0}=u_0\in \dot H^1,\quad \partial_tu_{\restriction t=0}=u_1\in L^2,
\end{gathered}\right.
\end{equation}

 In this
part we prove the following lemma:
\begin{lemma}
 \label{L:channel1}
Let $\beta>2$ and $u \in C(\R, \dot H^1) \cap L^5_t
L^{10}_x\left((-T, +T)\times \R^3\right)$ for any $T> 0$ be a
solution to equation \eqref{CP} with radial $(u_0, u_1)$. Suppose
for some $R>0$, we have
\begin{multline}
 \label{no_channel}
\lim_{t\to +\infty}\int_{|x|>R+|t|} |\nabla u(t,x)|^2+(\partial_tu(t,x))^2\,dx\\
=\lim_{t\to -\infty} \int_{|x|>R+|t|} |\nabla
u(t,x)|^2+(\partial_tu(t,x))^2\,dx=0.
\end{multline}
Then either $(u_0,u_1)$ or  $ (u_0,u_1)-\left(u_c,0\right)$ for some
$c\in \R$ is compactly supported.
\end{lemma}

 We first show two Lemmas.
\begin{lemma}
 \label{L:funct_decay} Let $u$ be as in Lemma \ref{L:channel1}, $v=ru$,
$v_0=ru_0$ and $v_1=ru_1$. If
\begin{equation}
 \label{small assu}
\int_{r_0}^{+\infty} \left((\partial_ru_0)^2+u_1^2\right)r^2\,dr\leq
\delta_0,
\end{equation}
for some large $r_0>0$  and small $\delta_0>0$, then
\begin{equation}
 \label{funct_decay}
\int_{r_0}^{+\infty} \Big((\partial_r v_0)^2+v_1^2\Big)\,dr\leq
C\frac{|v_0(r_0)|^{2}}{r^{2\beta-3}_0} + C
\frac{|v_0(r_0)|^{10}}{r^5_0}.
\end{equation}
Furthermore, for all $r,r'$ with $r_0\leq r\leq r'\leq 2r$,
\begin{equation}
 \label{bound}
\left|v_0(r)-v_0(r')\right|\leq C \frac{|v_0(r)|}{r^{\beta-2}} + C
\frac{|v_0(r)|^5}{r^2}\leq C \delta^2_0 |v_0(r)|.
\end{equation}
\end{lemma}

\smallskip
\noindent
{\bf Proof.}
We first prove \eqref{funct_decay}. Let $u^{L}(t,r)=S(t)(u_0,u_1)$
and $v^{L}=r u^{L}$. By Lemma \ref{lm:channelofenergyforlinearwave}, the following holds for all
$t\geq 0$ or for all $t\leq 0$
\begin{equation}
 \label{lin chan}
\int_{r_0+|t|}^{+\infty} (\partial_r v^{L}(t,r))^2+(\partial_t
v^{L}(t,r))^2\,dr\geq \frac 12 \int_{r_0}^{+\infty} (\partial_r
v_0(r))^2+v_1^2(r)\,dr.
\end{equation}

Let $\widetilde{u}$ be the solution of the following equation
\begin{equation*}
\left\{ \aligned
    \partial^2_{t} \widetilde{u }-  \Delta \widetilde{u} -\widetilde{V}(t,x) \widetilde{u} + \widetilde{u}^5 = &\; 0, \\
     (\widetilde{u}_0,\widetilde{u}_1)=\Psi_{r_0}(u_0,u_1)
\endaligned
\right.
\end{equation*}
and $\widetilde{u}^{L}=S(t)(\widetilde{u}_0,\widetilde{u}_1)$, where
$\widetilde{V}(t,x)$ is $ 0$ if $ |x|< |t|+r_0$, and $V(x)$ if $
|x|\geq |t|+r_0$. $\Psi_{r_0}(u_0,u_1)=(u_0,u_1)$ if $|x|\ge r_0$ and $=(u_0(r_0),0)$ if $|x|< r_0$. By \eqref{small assu}, we have
$\|(\widetilde{u}_0,\widetilde{u}_1)\|_{\dot H^1\times L^2}^2\leq
\delta_0|S^2|$. In addition, by the assumption on $V$, we know that $$
\big\|\widetilde{V}\big\|_{L^{\frac54}_tL^{\frac52}_x\left(\mathbb{R}\times
\mathbb{R}^3\right)}\lesssim r^{-\beta+2}_0. $$ Thus taking
$\delta_0$ small, $r_0$ large, we can get by the small data theory that for
all $t\in \R$,
\begin{align*}
 \left\| (\widetilde{u}-\widetilde{u}^{L},\partial_t \widetilde{u}-\partial_t\widetilde{u}^{L})(t)\right\|_{\dot H^1\times L^2}
 \leq   \frac{C}{r^{\beta-2}_0} \|(\widetilde{u}_0,\widetilde{u}_1)\|_{\dot H^1\times L^2} + C \|(\widetilde{u}_0,\widetilde{u}_1)\|^5_{\dot H^1\times
 L^2}.
\end{align*}
Since
$$\|(\widetilde{u}_0,\widetilde{u}_1)\|^2_{\dot H^1\times L^2} = |S^2| \left(\int_{r_0}^{+\infty}\left((\partial_r
v_0)^2+v_1^2\right)\,dr+r_0u_0^2(r_0)\right), $$
we have
\begin{multline*}
 \int_{r_0+|t|}^{+\infty}
  \left(\left(\partial_r \widetilde{u}^{L}(t)\right)^2+\left(\partial_t \widetilde{u}^{L}(t)\right)^2\right)r^2\,dr
\\
\leq 3\int_{r_0+|t|}^{+\infty} \left(\left(\partial_r
\widetilde{u}(t)\right)^2+\left(\partial_t
\widetilde{u}(t)\right)^2\right)r^2\,dr+
 3\frac{C}{r^{2\beta-4}_0}\left(\int_{r_0}^{+\infty}\left((\partial_r
v_0)^2+v_1^2\right)\,dr+r_0u_0^2(r_0)\right) \\
+ 3 C \left(\int_{r_0}^{+\infty}\left((\partial_r
v_0)^2+v_1^2\right)\,dr+r_0u_0^2(r_0)\right)^5 .
\end{multline*}
By the finite speed of propagation,
$$ \overrightarrow{u}(t,r)=\overrightarrow{\widetilde{u}}(t,r)\text{ and } \overrightarrow{u}^{L}(t,r)=\overrightarrow{\widetilde{u}}^{L}(t,r),\quad r\geq r_0+|t|,$$
and we obtain:
\begin{multline}
\label{nice_bound}
  \int_{r_0+|t|}^{+\infty} \left(\left(\partial_r u^{L}(t)\right)^2+\left(\partial_t u^{L}(t)\right)^2\right)r^2\,dr
\\
\leq  3\int_{r_0+|t|}^{+\infty} \left(\left(\partial_r
u(t)\right)^2+\left(\partial_t u(t)\right)^2\right)r^2\,dr+
 3\frac{C}{r^{2\beta-4}_0}\left(\int_{r_0}^{+\infty}\left((\partial_r
v_0)^2+v_1^2\right)\,dr+r_0u_0^2(r_0)\right) \\
+ 3 C \left(\int_{r_0}^{+\infty}\left((\partial_r
v_0)^2+v_1^2\right)\,dr+r_0u_0^2(r_0)\right)^5 .
\end{multline}
Combining \eqref{lin chan} and \eqref{nice_bound}, we see that the
following holds for all $t\geq 0$ or for all $t\leq 0$:
\begin{multline}
\label{nice_bound2}
 \frac 12 \int_{r_0}^{+\infty} \left((\partial_r v_0)^2+v_1^2\right)\,dr \\
\leq  3\int_{r_0+|t|}^{+\infty} \left(\left(\partial_r
u(t)\right)^2+\left(\partial_t u(t)\right)^2\right)r^2\,dr+
 3\frac{C}{r^{2\beta-4}_0}\left(\int_{r_0}^{+\infty}\left((\partial_r
v_0)^2+v_1^2\right)\,dr+r_0u_0^2(r_0)\right) \\
+ 3 C \left(\int_{r_0}^{+\infty}\left((\partial_r
v_0)^2+v_1^2\right)\,dr+r_0u_0^2(r_0)\right)^5 .
\end{multline}
Letting $t\to +\infty$ or $t\to -\infty$ in \eqref{nice_bound2}, we
see that the first term of the right-hand side of
\eqref{nice_bound2} goes to $0$ by \eqref{no_channel}.  Noting that
\eqref{small assu} and $r_0 u^{2}_0(r_0)=\frac{v_0^{2}(r_0)}{r_0}$,
we get \eqref{funct_decay}.

 We next prove \eqref{bound}. If $r_0\leq r\leq r'\leq 2r$, we have by \eqref{funct_decay}:
\begin{multline*}
 |v_0(r)-v_0(r')|\leq
 \sqrt{r}\sqrt{\int_{r}^{+\infty}(\partial_r
 v_0(\sigma))^2\,d\sigma}
 \leq C\frac{|v_0(r)|}{r^{\beta-2}} + C\frac{|v_0(r)|^5}{r^2} \leq
 C \delta^2_0 |v_0(r)|
\end{multline*}
which yields \eqref{bound}.\\

\smallskip
\noindent
\begin{lemma}
\label{L:limit} The function $v_0(r)$ has a limit $\ell\in \R$ as
$r\to +\infty$. Furthermore, there exists $C>0$ such that
\begin{equation}
\label{limit}
 \forall r\geq 1,\quad |v_0(r)-\ell|\leq \frac{C}{r^{\beta-2}}+ \frac{C}{r^2}.
\end{equation}
\end{lemma}

\smallskip
\noindent
{\bf Proof.}
 We first claim that there exists $C>0$ and $0<\epsilon < \min\{\beta-2,\frac{1}{10}\}$ such that for large $r$
\begin{equation}
 \label{weak_bound_v}
|v_0(r)|\leq Cr^{\epsilon}.
\end{equation}
Indeed by \eqref{bound}, if $n\in \mathbb{N}$,
$|v_0(2^{n+1}r_0)|\leq 2^{\epsilon} |v_0(2^nr_0)|$ if we choose $r_0$ sufficiently large. Hence by an
elementary induction
$$ |v_0(2^nr_0)|\leq 2^{n\epsilon}|v_0(r_0)|,$$
which shows the inequality \eqref{weak_bound_v} for $r=2^nr_0$,
$n\in \mathbb{N}$. The general case for \eqref{weak_bound_v} follows
from \eqref{bound}.

We next prove that $v_0(r)$ has a limit as $r\to+\infty$. By
\eqref{bound}, we get, for $n\in \mathbb{N}$,
$$ \left|v_0(2^nr_0)-v_0(2^{n+1}r_0)\right|\leq C \frac{\left|v_0(2^nr_0)\right|}{(2^nr_0)^{\beta-2}}+C\frac{\left|v_0(2^nr_0)\right|^5}{\left(2^nr_0\right)^2}.$$
By \eqref{weak_bound_v}, there exists $C>0$ such that
$$ \left|v_0(2^nr_0)-v_0(2^{n+1}r_0)\right|\leq C \frac{(2^nr_0)^{\epsilon}}{(2^nr_0)^{\beta-2}} + C\frac{(2^nr_0)^{5\epsilon}}{(2^nr_0)^{2}}.$$
Using that $\sum  \left( 2^{-(\beta-2-\epsilon)n} +
2^{-(2-5\epsilon)n} \right)$ converges for $0< \epsilon <\beta-2$,
we get
$$ \sum_{n\geq 1} \left|v_0(2^nr_0)-v_0(2^{n+1}r_0)\right|<\infty,$$
which shows that there exists $\ell \in \R$ such that
$$ \lim_{n\to+\infty} v_0(2^nr_0)=\ell.$$
Using \eqref{bound} and  \eqref{weak_bound_v}, we get $
\lim_{r\to+\infty} v_0(r)=\ell.$

It remains to prove \eqref{limit}. Using that $v_0(r)$ converges as
$r\to \infty$, we get that it is bounded for $r\geq r_0$, and thus
the first inequality in \eqref{bound} implies, for $r\geq r_0$ and
$n\in \mathbb{N}$,
$$\left|v_0(2^{n+1}r)-v_0(2^nr)\right|\leq \frac{C}{2^{(\beta-2)n}r^{\beta-2}} + \frac{C}{2^{2n}r^2}.$$
Summing up, we get
$$ \left|\ell - v_0(r)\right|=\left|\sum_{n\geq 0}(v(2^{n+1}r)-v(2^nr))\right|\leq \frac{C}{r^{\beta-2}}+ \frac{C}{r^2} ,$$
which concludes the proof of Lemma \ref{L:limit}.\\

\smallskip
\noindent
{\bf Proof of Lemma \ref{L:channel1}}.
 Consider the limit $\ell$ of $v_0$ defined in Lemma \ref{L:limit}. We distinguish between two cases, depending on $\ell$.

\noindent{\it Case: $\ell=0$} In this case we will show that
$(v_0,v_1)$ is compactly supported. We fix a large $r$. By
\eqref{bound}, if we take $r$ sufficiently large we have
$$\left|v_0(2^{n+1}r)\right|\geq \max(2^{-\frac{\beta-2}{2}}, \frac34)\left|v_0(2^nr)\right|,\quad \forall n\in \mathbb{N}.$$
By induction, we obtain $|v_0(2^nr)|\geq
\max(2^{-\frac{\beta-2}{2}n}, 3^n/4^n) |v_0(r)|$. Since $\ell=0$,
\eqref{limit} in Lemma \ref{L:limit} implies:
\begin{align*}
|v_0(2^nr)|\leq \frac{C}{(2^nr)^{\beta-2}} +  \frac{C}{4^nr^2}
\end{align*} Letting $n\to +\infty$, we get a contradiction unless
$v_0(r)=0$. Since $r$ is any large positive number, we have shown
that the support of $v_0$ is compact. By \eqref{funct_decay}, we get
that the support of $v_1$ is also compact, concluding this case.

\noindent{\it Case: $\ell\neq 0$} In this case we will show that
there exists a steady state solution $(u_c,0)$ such that
$\left(u_0- u_c,u_1\right)$ is compactly supported. We note that
for large $r$,
$$\left|u_c(r)-\frac{c}{r}\right|= o \left(\frac{1}{r}\right).$$
Thus Lemma \ref{L:limit} implies \begin{align} \label{u0 close to
uc}\left|u_c(r)-u_0(r)\right|\leq \left|u_c(r)-\frac{c}{r} \right| +
\left| u_0(r)-\frac{c}{r} \right|\leq
o\left(\frac{1}{r}\right),\quad r\geq 1,
\end{align} where $c=l$. Let $h=u-u_c$, $H=rh$.  Take a large
positive number $R_0$. Let $r_0>R_0$, and consider the solution $g$
of
\begin{equation}
\label{CP_g3} \left\{
\begin{gathered}
\partial_t^2 g -\Delta g-\widetilde{V}g+ 5U^4g+10U^3g^2+10U^2g^3+5Ug^4+g^5=0,\quad (t,x)\in \R\times \R^3\\
(g,\partial_t g)_{\restriction t=0}=(g_0,g_1)=\Psi_{r_0}(h_0,h_1).
\end{gathered}\right.
\end{equation}
where $U(t,x)$ is $0$ if $|x|< |t|+r_0$, and $u_c(x)$ if $|x|\geq
|t|+r_0$. By the finite speed of propagation,  we get
$$ \vec{g}(t,r)=\vec{h}(t,r)\text{ for } r\geq r_0+|t|.$$
Denote $\beta_1=\min\{\beta,4\}$. By \eqref{no_channel} and the similar proof of Lemma
\ref{L:funct_decay},  we can obtain that for a large $R_0>0$
\begin{align}
 \label{H funct decay}
\forall\; r_0>R_0,\quad \int_{r_0}^{+\infty}
\left((\partial_rH_0)^2+H_1^2\right)\,dr\leq
C\frac{H_0(r_0)^2}{r^{2\beta_1-3}_0}+ C \frac{H_0(r_0)^{10}}{r^5_0},
\end{align} which together with the assumption \eqref{small assu}
implies that
\begin{align} \big|H_0(r) - H_0(r')\big| \leq
C\frac{\big|H_0(r)\big|}{r^{\beta_1-2}}+ C \frac{H_0(r)^{5}}{r^2}\leq
C \delta^2_0 |H_0(r)|,
\end{align}
where $r_0 \leq r < r' < 2r$ and  $(H_0,H_1)=(H,\partial_t
H)_{\restriction t=0}$. Also $\lim_{r\to\infty}H_0(r)=0$. Analogous to the case $l=0$, we can obtain
$(H_0, H_1)$ is compactly supported.\\
\\

\bigskip
\noindent
Hao Jia and Baoping Liu, Department of Mathematics, University of Chicago.\\
Guixiang Xu, (visitor to University of Chicago) Institute of Applied Physics and Computational Mathematics, China.\\

\bigskip
\noindent
{\bf Acknowledgement.}\\
G. Xu was partly supported by the NSF of China (No. 11171033, No. 11231006) and was supported by the State Scholarship Fund, China Scholarship Council.\\

\smallskip
\noindent
 We are grateful to Prof. Schlag for suggesting the problem and to Prof. Kenig for encouragement. We also thank Prof. Polacik for telling us the ``transversality theorem" used in the proof of Theorem \ref{th:numberofsteadystates}, and Tianling Jin for the uniqueness of positive ground state in our situation.

\end{section}

\end{document}